\documentclass[a4paper,10pt,draft]{article}

\usepackage{amssymb, amsmath, amsthm, latexsym, verbatim, enumerate}

\newcommand \T{\theta}
\renewcommand \a{\alpha}
\newcommand \e{\epsilon}
\newcommand \K{\delta}

\newcommand \la{\lambda}
\newcommand \id{\mathrm{id}}

\newcommand \rk{\mathrm{rk}}
\newcommand \Ker{\mathrm{Ker \,}}

\newcommand \Span{\mathrm{Span}}
\newcommand \Tr{\mathrm{Tr} \, }

\newcommand \diag{\mathrm{diag}}

\newcommand \<{\langle}
\renewcommand \>{\rangle}

\newcommand \g{\mathfrak{g}}
\newcommand \h{\mathfrak{h}}
\newcommand \n{\mathfrak{n}}
\newcommand \m{\mathfrak{m}}
\newcommand \z{\mathfrak{z}}
\renewcommand \b{\mathfrak{b}}
\newcommand \ag{\mathfrak{a}}
\newcommand \ad{\mathrm{ad}}
\newcommand \ric{\mathrm{ric}}
\newcommand \Der{\mathrm{Der}}
\newcommand \End{\mathrm{End}}
\newcommand \F{\mathfrak{F}}
\newcommand \og{\mathfrak{o}}
\newcommand \p{\mathfrak{p}}
\newcommand \f{\mathfrak{f}}
\renewcommand \t{\mathfrak{t}}
\newcommand \s{\mathfrak{s}}
\renewcommand \l{\mathfrak{l}}
\newcommand \Og{\mathfrak{O}}
\renewcommand \d{\mathfrak{d}}

\theoremstyle{plane}
\newtheorem{theorem}{Theorem}
\newtheorem*{theorem*}{Theorem}
\newtheorem*{corollary*}{Corollary}
\newtheorem{lemma}{Lemma}
\newtheorem{proposition}{Proposition}
\newtheorem*{proposition*}{Proposition}

\newtheorem*{namedtheorem}{\theoremname}
\newcommand{\theoremname}{te}

\theoremstyle{definition}
\newtheorem{definition}{Definition}

\theoremstyle{remark}
\newtheorem{remark}{Remark}
\newtheorem{example}{Example}

\begin{document}

\title{Nilradicals of Einstein solvmanifolds}

\author{Y.Nikolayevsky}

\date{\today}

\maketitle

\begin{abstract}
A Riemannian Einstein solvmanifold is called standard, if the orthogonal complement to the nilradical of its Lie
algebra is abelian. No examples of nonstandard solvmanifolds are known.
We show that the standardness of an Einstein metric solvable Lie algebra is completely detected by its nilradical
and prove that many classes of nilpotent Lie algebras (Einstein nilradicals, algebras with less than four generators,
free Lie algebras, some classes of two-step nilpotent ones) contain no nilradicals of
nonstandard Einstein metric solvable Lie algebras. We also prove that there are no nonstandard Einstein metric
solvable Lie algebras of dimension less than ten.
\end{abstract}

\section{Introduction}
\label{s:intro}

The theory of Riemannian homogeneous spaces with an Einstein metric splits into three very different cases
depending on the sign of the Einstein constant, the scalar curvature. Among them, only in the Ricci-flat case
the picture is complete: by the result of \cite{AK}, every Ricci-flat homogeneous space is flat.
Despite many remarkable existence, nonexistence and partial classification results (see \cite{WZ, BWZ}), the positive
scalar curvature case is far from being completely understood.
One of the main open questions in the case of negative scalar curvature (which we are dealing with in this paper),
is the \emph{Alekseevski Conjecture} \cite{A1} asserting that any Einstein homogeneous Riemannian space
with negative scalar curvature admits a simply transitive solvable isometry group. This is equivalent to saying
that any such space is a \emph{solvmanifold}, a solvable Lie group with a left-invariant Riemannian metric satisfying
the Einstein condition.

Another question arises from the fact that all the known examples of Einstein solvmanifolds are \emph{standard}.
This means that the metric solvable Lie algebra $\g$ of such a solvmanifold has the following property: the orthogonal
complement $\ag$ to the derived algebra of $\g$ is abelian.
\emph{Are there any nonstandard Einstein solvmanifolds?}

In the paper \cite{H}, J. Heber gave a deep and detailed analysis of standard Einstein solvmanifolds and also
showed that for several classes of solvmanifolds, the existence of an Einstein left invariant
metric implies standardness. What is more, nonstandard solvmanifolds must be algebraically very different from the
standard ones. In particular, a solvable Lie group cannot carry two Einstein left-invariant metrics, one of which is
standard, and the other one is not. Moreover, standard Einstein solvmanifolds form an open and compact subspace in
the moduli space of all Einstein solvmanifolds of the given dimension and fixed scalar curvature.

Further progress was made by D.Schueth in \cite{S}, for solvmanifolds with $\dim[\ag, \ag] = 1$. All the Einstein
solvmanifolds of dimension up to $7$ are standard, as it is shown by Yu.Nikpnorov: for dimension up to $6$, this
follows from their classification obtained in \cite{N1, NN}, for Einstein solvmanifolds of dimension $7$, the
standardness is proved in \cite{N2}.

The problem of studying and classifying standard Einstein solvmanifolds (and standard Einstein metric solvable Lie
algebras) was given a great deal of attention in the last decade
(a very incomplete list of references is \cite{GK, L1, L2, L3, LW, Ni, N1, P, T}).

This paper addresses the following question: \emph{which algebraic properties of a nilpotent Lie
algebra $\n$ force any Einstein metric solvable Lie algebra with the nilradical $\n$ to be standard?}
The focus on the nilradical is motivated by the following two reasons. Firstly, as it is shown by J.Lauret \cite{L1},
in the standard case, the nilradical is the only
thing one needs to know to reconstruct both the Lie algebra structure and the inner product of any of its Einstein
metric solvable extension (see Section~\ref{s:facts} for details). Secondly, if the nonstandard Einstein
solvmanifolds do exist at all, a natural approach to finding one might be to start with a nilpotent Lie algebra,
which will then serve as the nilradical for the Lie algebra of such a solvmanifold.

The idea of using the properties of the nilradical to study the standardness also appeared
in the paper \cite{N2}, one of the main result of which is that if for a metric nilpotent Lie
algebra, the symmetric operator $\id - (\Tr(\ric)/\Tr (\ric^2)) \, \ric$ is positive, then any
Einstein metric solvable extension of it is standard.

\bigskip

An endomorphism $A$ of a linear space is called \emph{semisimple}, if its nilpotent part vanishes, is called \emph{real},
if all its eigenvalues are real, and is called \emph{nonnegative} (respectively, \emph{positive}), if all its eigenvalues
are real and nonnegative (respectively, positive). In the latter case, we write $A \ge 0$ (respectively, $A>0$).
A derivation $\phi$ of a Lie algebra $\l$ is called \emph{pre-Einstein}, if it is real, semisimple, and
\begin{equation*}
    \Tr (\phi \circ \psi) = \Tr \psi,  \quad \text{for any $\psi \in \Der(\l)$}.
\end{equation*}
In Proposition~\ref{p:pEexistunique} we show that \emph{every Lie algebra $\l$ admits a unique pre-Einstein derivation
$\phi_\l$, up to a conjugation} (it could well happen that $\phi_\l = 0$, even for a nilpotent $\l$).

For any $\psi \in \Der(\l)$ we denote $\ad_\psi$ the corresponding inner derivation of $\Der(\l)$. If $\psi$ is
semisimple and real, the same is true for $\ad_\psi$.

The main technical result of the paper is the following sufficient condition for the standardness.

\begin{theorem}\label{th:pE}
Let $\n$ be a nilpotent Lie algebra with $\phi$ a pre-Einstein derivation. Suppose that
\begin{equation}\label{eq:pE>0}
  \phi > 0 \qquad \text{and} \qquad \ad_\phi \ge 0.
\end{equation}
Then any Einstein metric solvable Lie algebra with the nilradical $\n$ is standard.
\end{theorem}

This theorem (and most of the subsequent ones) can be slightly strengthened: given an Einstein metric solvable Lie
algebra $(\g,\<\cdot,\cdot\>)$, it suffices to choose $\n$ anywhere between its derived algebra and its nilradical
and to check the condition \eqref{eq:pE>0} (see the proof of Theorem~\ref{th:pE} in Section~\ref{s:thpE}).

It would be tempting to get rid of the second condition ($\ad_\phi \ge 0$) in Theorem~\ref{th:pE}. Although we do not
know how to avoid it completely, the following can be proved:

\begin{theorem}\label{th:split}
Let $\n$ be a nilpotent Lie algebra whose pre-Einstein derivation is positive.
Then any \emph{split} Einstein metric solvable Lie algebra with the nilradical $\n$ is standard.
\end{theorem}

Recall that a Lie algebra $\mathfrak{l}$ is called \emph{split}, if for any $X \in \mathfrak{l}$ there exists
$Y \in \mathfrak{l}$ such that $\ad_Y$ is the nilpotent part of $\ad_X$. For example, any algebraic Lie algebra is split.

As the main consequence of Theorem~\ref{th:pE}, we show that for an Einstein metric solvable Lie algebra, it is the
nilradical which detects the standardness:

\begin{theorem} \label{th:nil}
Let $\g$ and $\g_0$ be solvable Lie algebras with the same nilradical. If $\g_0$ admits a standard Einstein inner
product, then any Einstein inner product on $\g$ is also standard.
\end{theorem}

The idea of this theorem is close to \cite[Theorem 5.3]{H} saying that a solvable Lie algebra cannot carry two
Einstein inner products, one of which is standard, and another one is not. Theorem~\ref{th:nil} gives the positive
answer to the question asked in \cite{N2}.

A nilpotent Lie algebra which can be a nilradical of a \emph{standard} Einstein metric solvable Lie algebra is called
an \emph{Einstein nilradical}.
If $\n$ is an Einstein nilradical, its pre-Einstein derivation is a positive multiple of the Einstein derivation
(see the proof of Theorem~\ref{th:nil}).
Not every nilpotent Lie algebra is an Einstein nilradical. By \cite{H}, a necessary condition for that is to admit
an $\mathbb{N}$-gradation (which is defined by the Einstein derivation). This condition is, however, far from being
sufficient, and a very delicate analysis (see the stratification constructed in \cite{LW}) is required to decide
whether a given $\mathbb{N}$-graded nilpotent Lie algebra is indeed
an Einstein nilradical (examples of those which are not can found in \cite{LW, Ni, P}).

Each of the following nilpotent Lie algebras $\n$ is an Einstein nilradical (see \cite{L4} for a complete
up-to-date list), hence by Theorem~\ref{th:nil}, any Einstein metric solvable Lie algebra with such a nilradical $\n$
is standard:

\begin{enumerate}[(i)]
  \item $\n$ is abelian \cite{A2};
  \item $\n$ has a codimension one abelian ideal \cite{L3};
  \item $\dim \n \le 6$; \cite{W, L3};
  \item $\n$ is uniform two-step nilpotent (this class includes, for example, nilpotent Lie algebras of Heisenberg
  type, that is, those arising from representations of Clifford algebras) \cite{GK};
  \item $\n$ is a direct sum of Einstein nilradicals \cite{P};
  \item $\n = \t_n$, the algebra of strictly upper-triangular $n \times n$ matrices \cite{P}.
\end{enumerate}

The fact that any Einstein metric solvable Lie algebra with an abelian nilradical $\n$ is standard is independently
proved in \cite{N2}.

In view of (i), one may wonder what happens in the other extremal case, when $\n$ is the least commutative, a filiform.
The answer is given by the following theorem (recall that a filiform algebra is generated by two elements):

\begin{theorem}\label{th:filiform}
Let $\n$ be a nilpotent Lie algebra generated by no more than three elements (that is, $\dim \n - \dim [\n, \n] \le 3$).
Then any Einstein metric solvable Lie algebra with the nilradical $\n$ is standard.
\end{theorem}

Another consequence of Theorem~\ref{th:pE} is the following theorem:

\begin{theorem}\label{th:free}
Any Einstein metric solvable Lie algebra with a free nilradical is standard.
\end{theorem}

Free Einstein nilradical were classified in \cite{Ni}. It is shown that apart from the abelian and the two-step ones,
only six other free Lie algebras can be Einstein nilradicals.

In the second part of the paper, we consider the applications of Theorem~\ref{th:pE} and Theorem~\ref{th:nil}
to two-step nilpotent nilradical. A two-step nilpotent Lie algebra $\n$ is said to
be of type $(p, q)$, if $\dim \n = p+q$ and $\dim [\n, \n] = p$ (clearly, $p\le D:=q(q-1)/2$). Any such algebra is
determined by a $p$-dimensional subspace $W$ of $q \times q$ skew-symmetric
matrices, that is, by a point on the Grassmannian $G(p, \og(q))$.
Two algebras of type $(p, q)$ are isomorphic if and only if they lie in the same orbit $[W]$ of the action of
the group $\mathrm{SL}(q)$ on  $G(p, \og(q))$.

Let $\Og(p, q)$ be the class of two-step nilpotent Lie algebras $\n$ of
type $(p, q)$, which have the following property:
for any $\psi \in \Der(\n),\quad \Tr \psi = 0 \Rightarrow \Tr \psi_{|[\n,\n]} = 0$ (an equivalent definition: the
canonical derivation acting as multiplication by $2$ on $[\n, \n]$ and as the identity on some linear complement to
$[\n, \n]$ is pre-Einstein).
The class $\Og(p, q)$ consist of the two-step nilpotent algebras whose derivation algebra has, in the sense, the nicest
possible structure. Another natural motivation to single the class $\Og(p, q)$ out
is that the Einstein nilradicals belonging to $\Og(p, q)$ are
precisely those whose eigenvalue type is $(1 < 2; q, p)$ which is, perhaps, the most well-studied class of the
Einstein nilradicals. We prove the following ``general position" theorem:

\begin{theorem}\label{th:twostepopen}
1. No two-step nilpotent Lie algebra $\n \in \Og(p, q)$ can be the nilradical of a nonstandard Einstein metric
solvable Lie algebra.
For all the pairs $(p, q)$ such that $1 < p < D-1$, except when $q$ is odd and $p=2$ or $p=D-2$, the set
$\{W \,|\, [W] \in \Og(p,q)\}$ contains an open and dense subset of $G(p, \og(q))$.
The set $\Og = \cup_{p,q}\Og(p, q)$ is closed under taking direct products and passing to the
dual algebra.

2. An Einstein metric solvable Lie algebra with a two-step nilradical $\n$ is standard in each of the following cases:
\begin{enumerate}[(a)]
  \item $\n$ is nonsingular;
  \item $\n$ is of type $(p,q), \; q \le 4$;
  \item $\n$ is of type $(D-1, q)$ (in fact, such an $\n$ is even an Einstein nilradical).
\end{enumerate}

\end{theorem}

Recall that a two-step nilpotent Lie algebra $\n$ is called \emph{nonsingular}, if for every
$X \in \n \setminus [\n, \n]$, the map $\ad_{X}: \n \to [\n, \n]$ is surjective.
As to the missing pairs
$(p, q)$ in assertion 1, note that any two-step nilpotent algebra with $p=1$ or $p=D$ is an Einstein nilradical.
The case $p = D-1$ is covered by (c) of assertion 2 (see Remark~\ref{r:missing} in Section~\ref{ss:2stepopq}).
As for almost all the pairs $(p, q)$, the action of $\mathrm{GL}(q)$ on
$G(p, \og(q))$ has no open orbits \cite{E1}, the set $\Og(p,q)$ is of the second Baire category 
in the non-Hausdorff space of orbits $G(p, \og(q))/\mathrm{SL}(q)$, the space of isomorphism classes of algebras of
type $(p, q)$.

Using Theorems~\ref{th:pE}, \ref{th:nil}, \ref{th:filiform} and \ref{th:twostepopen} we prove the following theorem:

\begin{theorem}\label{th:dim7}
Any Einstein metric solvable Lie algebra whose nilradical has dimension at most $7$ is standard. In particular,
any Einstein metric solvable Lie algebra of dimension less than ten is standard.
\end{theorem}

The paper is organized as follows. In Section~\ref{s:facts} we provide the necessary background on Einstein
solvmanifolds. In Section~\ref{s:thpE} we establish the existence and uniqueness of the pre-Einstein derivation
(Proposition~\ref{p:pEexistunique}) and give the proof of Theorems~\ref{th:pE} and \ref{th:split}. In
Section~\ref{s:appl} we consider the applications of Theorem~\ref{th:pE} to Einstein nilradicals
and to free nilradicals proving Theorems~\ref{th:nil} and \ref{th:free}, respectively. The proof
of Theorem~\ref{th:filiform} is given in Section~\ref{s:fili}. In Section~\ref{s:2step} we consider the two-step
nilradicals and prove Theorem~\ref{th:twostepopen}. Finally, in Section~\ref{s:dim7} we establish
Theorem~\ref{th:dim7}.

The author is grateful to Yu.G.Nikonorov for useful discussions.

\section{Preliminaries}
\label{s:facts}

Let $\g$ be a Lie algebra, with $B$ the Killing form. The coboundary operator $\delta$ acts from the space of
endomorphisms of (the linear space) $\g$ to the space of $\g$-valued two-forms on $\g$ by the formula
$\delta(A) (X, Y) := -A[X, Y] + [AX, Y] + [X, AY]$, where $A \in \End(\g),\; X, Y \in \g$. From the definition,
$\delta(\id) (X, Y) = [X, Y]$ and $\Ker \delta = \Der(\g)$, the algebra of derivations of $\g$.

For an inner product $\< \cdot, \cdot \>$ on $\g$, define the \emph{mean curvature vector} $H$ by
$\<H, X\> = \Tr \ad_X$ (clearly, $H$ is orthogonal to the nilradical of $\g$).
For $A \in \End(\g)$, let $A^*$ be its metric adjoint and $S(A) :=\frac12 (A +A^*)$ be the
symmetric part of $A$. Define a symmetric form $\Xi$ on $\End(\g)$ by
\begin{equation}\label{eq:Phidelta}
\Xi(A_1, A_2)= \<\delta(A_1), \delta(A_2)\> = \sum\nolimits_{i,j} \<\delta(A_1)(E_i, E_j), \delta(A_2)(E_i, E_j)\>,
\end{equation}
for an orthonormal basis $\{E_i\}$ for $\g$.
Then the Ricci operator $\ric$ of the metric Lie algebra $(\g, \< \cdot, \cdot \>)$
(the symmetric operator associated to the Ricci tensor)
is implicitly defined by
\begin{equation}\label{eq:riccidefdelta}
\Tr \bigl(\ric + S(\ad_H) + \frac12 B \bigr) \circ A
= -\frac14 \, \Xi(A, \id),
\end{equation}
for any $A \in \End(\g)$, or, in the expanded form, by
\begin{equation}\label{eq:riccidef}
\Tr \bigl(\ric + \frac12 (\ad_H + \ad_H^*) + \frac12 B \bigr) \circ A
= \frac14 \sum\nolimits_{i,j} \<A[E_i, E_j] - [AE_i, E_j] - [E_i, AE_j], [E_i, E_j]\>.
\end{equation}
For an Einstein metric solvable Lie algebra of negative scalar curvature $c \dim \g$ this gives
\begin{equation}\label{eq:riccidefein}
\Tr \bigl(|c| \, \id_\g - S(\ad_H) - \frac12 B \bigr) \circ A
= -\frac14 \sum\nolimits_{i,j} \<A[E_i, E_j] - [AE_i, E_j] - [E_i, AE_j], [E_i, E_j]\>.
\end{equation}

If $(\n, \< \cdot, \cdot \>)$ is a nilpotent metric Lie algebra, then $H = 0$ and $B = 0$, so \eqref{eq:riccidef} gives
\begin{equation}\label{eq:riccinil}
\Tr (\ric_{\n} \circ A) = \frac14 \, \Xi(A, \id)
= \frac14 \sum\nolimits_{i,j} \<A[E_i, E_j] - [AE_i, E_j] - [E_i, AE_j], [E_i, E_j]\>,
\end{equation}
for any $A \in \End(\n)$.

\begin{definition} \cite{H}
\label{d:st}
An inner product on a solvable Lie algebra $\g$ is called \emph{standard}, if the orthogonal complement to
the derived algebra $[\g, \g]$ is abelian. A metric solvable Lie algebra $(\g, \<\cdot,\cdot\>)$ is called
\emph{standard}, if the inner product $\<\cdot,\cdot\>$ is standard.
\end{definition}

As it is proved in \cite{AK}, any Ricci-flat metric solvable Lie algebra is flat. By the result of \cite{DM},
any Einstein metric solvable unimodular Lie algebra is also flat. In what follows, we always assume $\g$
to be nonunimodular ($H \ne 0$), with an inner product of a strictly negative scalar curvature $c \dim \g$.

Any standard Einstein metric solvable Lie algebra admits a rank-one reduction \cite[Theorem 4.18]{H}. This means that if
$(\g, \< \cdot, \cdot\>)$ is such an algebra, with the nilradical $\n$ and the mean curvature vector $H$, then the
subalgebra $\g_1 = \mathbb{R}H \oplus \n$, with the induced inner product, is also Einstein and standard. What is
more, the derivation $\Phi=\ad_{H|\n}:\n \to \n$ is symmetric with respect to the inner product, and all its
eigenvalues belong to $\a \mathbb{N}$ for some constant $\a > 0$. This implies, in particular, that the nilradical $\n$
of a standard Einstein metric solvable Lie algebra admits an $\mathbb{N}$-gradation defined by the eigenspaces of $\Phi$.
As it is proved in \cite[Theorem 3.7]{L1}, a necessary and sufficient condition for a metric nilpotent algebra
$(\n, \< \cdot, \cdot\>)$ to be the nilradical of a standard Einstein metric solvable Lie algebra is
\begin{equation}\label{eq:ricn}
    \ric_\n = c \, \id_\n + \Phi,
\end{equation}
where $c \dim \g < 0$ is the scalar curvature of $(\g, \< \cdot, \cdot\>)$. This equation, in fact, defines
$(\g, \< \cdot, \cdot\>)$ in the following sense: given a metric nilpotent Lie algebra whose Ricci operator
satisfies \eqref{eq:ricn}, with some constant $c < 0$ and some $\Phi \in \Der(\n)$, one can define $\g$ as a
one-dimensional extension of $\n$ by $\Phi$. For such an extension $\g = \mathbb{R}H \oplus \n, \; \ad_{H|\n} = \Phi$,
and the inner product
defined by $\<H, \n \> = 0,\; \|H\|^2 = \Tr \Phi$ (and coinciding with the existing one on $\n$) is Einstein, with the
scalar curvature $c \dim \g$. Following \cite{L1} we call a nilpotent Lie algebra $\n$ which admits an inner product
$\< \cdot, \cdot\>$ and a derivation $\Phi$ satisfying \eqref{eq:ricn} an \emph{Einstein nilradical}, the
corresponding derivation $\Phi$ is called an \emph{Einstein derivation}, and the inner product $\< \cdot, \cdot\>$
the \emph{nilsoliton metric}.

As it is proved in \cite[Theorem 3.5]{L1}, a nilpotent Lie algebra admits no more than one nilsoliton metric, up to
conjugation and scaling (and hence, an Einstein derivation, if it exists, is unique, up to conjugation and scaling).
Equation \eqref{eq:ricn}, together with \eqref{eq:riccinil}, implies that if $\n$ is an Einstein nilradical, with
$\Phi$ the Einstein derivation, then
\begin{equation}\label{eq:tracestandard}
    \Tr (\Phi \circ \psi) = - c \, \Tr \psi,  \quad \text{for any $\psi \in \Der(\n)$}.
\end{equation}

We use the following notational convention: for a linear space $V, \quad \End(V)$ is the Lie algebra of linear
transformations of $V$ (that is, $\End(V) = \mathfrak{gl}(V)$), $\oplus$ is the direct sum of linear spaces (even
if the summands are Lie algebras). For an endomorphism $A$ of a linear space $V$, we denote $A^S$ and $A^N$ its
semisimple and nilpotent parts respectively: $A = A^S + A^N, \; [A^S, A^N] = 0$. $A^S$ can be further decomposed as
$A^S=A^{\mathbb{R}} + A^{i\mathbb{R}}$, the real and the imaginary part of $A^S$, respectively. The operator
$A^{\mathbb{R}}$ is defined as follows: if $V_1, \ldots, V_m$ are the eigenspaces of $A^S$ acting on $V^\mathbb{C}$,
with the eigenvalues $\la_1, \ldots, \la_m$, respectively, then $A^{\mathbb{R}}$ acts by multiplication by
$\mu \in \mathbb{R}$ on every subspace $(\oplus_{k: \mathrm{Re} \la_k = \mu} V_k) \cap V$. For any $A \in \End(V)$, the
operators $A, A^S, A^N, A^{\mathbb{R}}$, and $A^{i\mathbb{R}}$ commute. If $\psi$ is a derivation of a Lie algebra
$\g$, then each of the $\psi^S, \psi^N, \psi^{\mathbb{R}}$, and $\psi^{i\mathbb{R}}$ is also a derivation.

\section{Pre-Einstein Derivation. Proof of Theorems~\ref{th:pE} and \ref{th:split}}
\label{s:thpE}

In this section we prove Theorem~\ref{th:pE} and Theorem~\ref{th:split} and consider some examples.

\begin{definition} \label{d:pE}
A derivation $\phi$ of a Lie algebra $\l$ is called \emph{pre-Einstein}, if it is real,
semisimple, and
\begin{equation}\label{eq:pEtrace}
    \Tr (\phi \circ \psi) = \Tr \psi,  \quad \text{for any $\psi \in \Der(\l)$}.
\end{equation}
\end{definition}

We start with the following proposition.

\begin{proposition} \label{p:pEexistunique}
{\ }

(a) An arbitrary Lie algebra $\l$ admits a pre-Einstein derivation $\phi_\l$.

(b) The derivation $\phi_\l$ is determined uniquely up to an automorphism of $\l$.

(c) All the eigenvalues of $\phi_\l$ are rational numbers.
\end{proposition}
\begin{proof}
(a) The algebra $\Der (\l)$ is algebraic. Let $\Der (\l) = \s \oplus \t \oplus \n$ be its Mal'cev decomposition, where
$\t \oplus \n$ is the radical of  $\Der (\l)$, $\s$ is semisimple, $\n$ is the set of all nilpotent elements in
$\t \oplus \n$ (and is the nilradical of $\t \oplus \n$), $\t$ is a torus, an abelian subalgebra consisting of
semisimple elements, and $[\t, \s] = 0$. With any $\psi \in \t,\; \psi^{\mathbb{R}}$ and $\psi^{i\mathbb{R}}$ are also
in $\t$. The subspaces
$\t_c=\{\psi^{\mathbb{R}}\,|\, \psi \in \t\}$ and $\t_s=\{\psi^{i\mathbb{R}}\,|\, \psi \in \t\}$ are the compact
and the completely $\mathbb{R}$-reducible tori (the elements of $\t_s$ are diagonal matrices in some basis for
$\l$), $\t_s \oplus \t_c = \t$.

The quadratic form $b$ defined on $\Der (\l)$ by $b(\psi_1, \psi_2) = \Tr (\psi_1 \circ \psi_2)$ is invariant
($b(\psi_1, [\psi_2, \psi_3]) = b([\psi_1, \psi_3], \psi_2))$. In general, $b$ is degenerate, with
$\Ker b = \n$, 
so for any $\psi \in \n,\; \b(\t, \psi) = \Tr \psi = 0$. As $\s$ is semisimple and $[\t, \s] = 0$, we also have
$b(\t, \psi) = \Tr \psi = 0$, for any $\psi \in \s$. Moreover, for any
$\psi \in \t_c, \quad b(\t_s, \psi) = \Tr \psi = 0$.

So to produce a pre-Einstein derivation for $\l$ it suffices to find an element $\phi \in \t_s$ which satisfies
\eqref{eq:pEtrace}, with any $\psi \in \t_s$. Such a $\phi$ indeed exists, as the restriction of $b$ to $\t_s$ is
nondegenerate (even definite) and is unique, when a particular torus $\t$ is chosen.

(b) The subalgebra $\s \oplus \t$ is a maximal reducible subalgebra of $\Der(\l)$. As by \cite[Theorem 4.1]{Mo},
the maximal reducible subalgebras of $\Der(\l)$ are conjugate by an inner automorphism of $\Der(\l)$ (which is an
automorphism of $\l$), and then $\t$, the center of $\s \oplus \t$, is defined uniquely, the uniqueness of $\phi$,
up to an automorphism, follows.

(c) The proof is similar to \cite[Theorem 4.14]{H}. Suppose $\phi$
has eigenvalues $\mu_i$, with multiplicities $d_i$, respectively, $i=1, \ldots, m$. In a Euclidean space $\mathbb{R}^m$
with a fixed orthonormal basis $f_i$, consider all the vectors of the form $f_i+f_j-f_k$ such that
$\mu_i+\mu_j-\mu_k = 0$.
In their linear span choose a basis $v_\a, \; \a = 1, \ldots , p$, consisting of vectors of the above form and
introduce an $m \times p$ matrix $F$ whose vector-columns are the $v_\a$'s. Then any vector
$\nu = (\nu_1, \ldots, \nu_m)^t \in \mathbb{R}^m$ satisfying $F^t\nu =0$ defines a derivation $\psi = \psi(\nu)$ having
the same eigenspaces as $\phi$, but with the corresponding eigenvalues $\nu_i$. From \eqref{eq:pEtrace} we must have
$\sum d_i (\mu_i -1) \nu_i = 0$ for any such $\nu$, which implies that the vector
$(d_1 (\mu_1 -1), \ldots, d_m (\mu_m -1))^t$ belongs to the column space of $F$. So there exists $x \in \mathbb{R}^p$
such that $\mu = 1_m + D^{-1}Fx$, where $\mu = (\mu_1, \ldots, \mu_m)^t$, $1_m = (1, \ldots, 1)^t \in \mathbb{R}^m$,
and $D = \diag (d_1, \ldots, d_m)$. As $\phi$ by itself is a derivation, we have $F^t \mu = 0$, which implies
$F^t1_m + F^tD^{-1}Fx = 0$, so that $x = - (F^tD^{-1}F)^{-1}1_p$, as $F^t1_m = 1_p$ and $\rk \, F = p$. Then
$\mu = 1_m - D^{-1}F(F^tD^{-1}F)^{-1}1_p$ and the claim follows, as all the entries of $D$ and of $F$ are integers.
\end{proof}

The following easy lemma shows that in the definition of the standardness, the derived algebra can be replaced by the
nilradical (in fact, by any ideal in between). Note that for a standard $\g$, the nilradical coincides with the
derived algebra \cite[Corollary 4.11]{H}.

\begin{lemma}\label{l:nilradvsderived}
Suppose that for an Einstein metric solvable Lie algebra $(\g, \<\cdot,\cdot\>)$, the orthogonal
complement $\ag$ to the nilradical $\n$ is abelian. Then $\g$ is standard.
\end{lemma}
\begin{proof}
The proof is close to that of \cite[Corollary 2.4]{S}.
By \cite[Lemma 4.11]{H}, for any $Y \in \ag$, the operator $\ad_Y^*$ is a derivation, which
must be normal by \cite[Lemma 2.1]{H}. As for any $X \in \n \cap [\g, \g]^\perp, \; \ad_Y^*X = 0$, the normality of
$\ad_Y$ implies $\ad_Y^*\ad_YX=0$, so $[Y, X] = 0$.
In particular, $[H, X] = 0$, so, again, $\ad_X^* \in \Der(\g)$ and $\ad_X$ is normal. As $\ad_X$ is nilpotent,
this implies $\ad_X = 0$, and so $X=0$ (as the Ricci curvature in the direction of $X$ is nonnegative, see e.g.
\cite[Lemma 2.1]{Mi}).
\end{proof}

The rest of the proof of Theorem~\ref{th:pE} follows the lines of the proof of \cite[Theorem 5.3]{H}:
given an Einstein metric solvable Lie algebra $(\g, \<\cdot, \cdot \>)$, we compare two derivations of $\n$:
the pre-Einstein derivation $\phi_\n$ and $\ad_{H|\n}$, and show that in the assumptions \eqref{eq:pE>0} they must
coincide (up to conjugation and scaling). Note that a similar approach was also used in the proof
of \cite[Theorem 4]{NN}.

We need the following technical lemma for references. The proof of assertion 1 is essentially contained in the proof
of \cite[Theorem 5.3]{H}, the proof of assertion 2 follows from \cite[Theorem 3]{N2} (compare \cite[Lemma 4.1]{H}).

\begin{lemma}\label{l:gs}
Let $(\l, \<\cdot, \cdot \>)$ be a metric Lie algebra and let $\psi \in \Der(\l)$.

1. Suppose $\psi$ is semisimple, real, with the eigenvalues $\la_1 < \ldots < \la_m$ and the corresponding eigenspaces
$\n_1, \ldots, \n_m$. Define a symmetric operator $\psi^{\perp}$ 
as follows: for $W_i=\oplus_{j=i}^m \n_j \cap (\oplus_{j=i+1}^m  \n_j)^{\perp}$, set
$\psi^{\perp}_{|W_i} = \la_i \, \id_{|W_i}$. Then $\Xi(\psi^\perp, \, \id) \le 0$,
with the equality if and only if $\psi^{\perp}$ is a derivation.

2. For an arbitrary $\psi \in \Der(\l)$, denote
$S = (\psi + \psi^*)/2$ the symmetric part of $\psi$. Then
\begin{equation}\label{eq:him}
    \Tr (\ric + \ad_H + \frac12 B) \circ [\psi, \psi^*]
    = \frac12 \, \Xi(S,S) = \frac12 \sum\nolimits_{ij} \|S[E_i, E_j] - [SE_i, E_j] - [E_i, SE_j]\|^2 \ge 0,
\end{equation}
where $\{E_i\}$ is an orthonormal basis for $(\l, \<\cdot, \cdot \>)$.
In particular, if $(\l, \<\cdot, \cdot \>)$ is an Einstein nilradical, with the Einstein derivation $\Phi$, then
\begin{equation}\label{eq:hinil}
\Tr (\Phi \circ [\psi, \psi^*]) = \frac12 \sum\nolimits_{ij} \|S[E_i, E_j] - [SE_i, E_j] - [E_i, SE_j]\|^2 \ge 0.
\end{equation}

3. Let $\phi_i$ be the pre-Einstein derivations for the Lie algebras $\l_i, \; i =1, 2$. Then the pre-Einstein derivation
$\phi_\l$ for the Lie algebra $\l$, the direct sum of the Lie algebras $\l_1$ and $\l_2$, is $\phi_1 \oplus \phi_2$.

\end{lemma}

\begin{proof}
1. As for any $i = 1, \ldots, m,\; \oplus_{s=i}^m W_s = \oplus_{s=i}^m \n_s$ and $[\n_i, \n_j] \subset \n_k$, with
$\la_k = \la_i + \la_j$ (or zero, if there is no such $k$), $[W_i, W_j] \subset \oplus_{s=k}^m W_s:= U_k$. The
subspace $U_k$ is invariant for $\psi^\perp$, and the restriction of $\psi^\perp$ to
it is a symmetric operator whose eigenvalues are greater than or equal to $\la_k = \la_i + \la_j$.
So for any $E_i \in W_i,\, E_j \in W_j, \quad
\<[\psi^\perp E_i, E_j] + [E_i, \psi^\perp E_j] - \psi^\perp[E_i, E_j], [E_i, E_j]\> =
\<((\la_i + \la_j)\id - \psi^\perp)_{|U_k} [E_i, E_j], [E_i, E_j]\> \le 0$,
with the equality only when $\psi^\perp$ is also a derivation.

2. The quadratic form $\Xi$ defined by \eqref{eq:Phidelta} on the linear space $\End(\l)$ is positive semidefinite,
with $\Ker \Xi = \mathrm{Der}(\l)$. Clearly, $\Xi(A, \, \id) = \Xi(A^*, \, \id)$. Moreover, for any symmetric $S$ and
a skew-symmetric $K$, we have
\begin{equation}\label{eq:sk}
\Xi(S, K) = -\Xi(SK, \,\id) = \Xi(KS, \,\id).
\end{equation}
To prove this, note that the expansion of $\Xi(S, K) = \<\delta(S), \delta(K)\>$ is a sum of nine terms,
six of which vanish, so that
\begin{equation*}
\Xi(S, K) = \sum\nolimits_{ij} (\<S[E_i,E_j], K[E_i, E_j]\> + \<[SE_i,E_j], [KE_i,E_j]\> + \<[E_i,SE_j], [E_i,KE_j]\>).
\end{equation*}
Indeed, for example,
$\sum_{ij} \<S[E_i,E_j], [KE_i, E_j]\>= \sum_{ij} \<S \ad_{E_j}E_i , \ad_{E_j}KE_i\> =
\sum_j \Tr \ad_{E_j}^* S \ad_{E_j} K = 0$, as the endomorphisms $\ad_{E_j}^* S \ad_{E_j}$ are self-adjoint. Also,
$\sum_{ij} \<[SE_i,E_j], [E_i, KE_j]\> = 0$, which can be seen if we take $E_i$ as the eigenvectors of $S$. The
vanishing of the other four terms can be shown by a similar routine check. Then \eqref{eq:sk} follows from
$\sum_{ij} \<[SE_i,E_j], [KE_i, E_j]\> = - \sum_j \Tr (K \ad_{E_j}^* \ad_{E_j} S)
= - \sum_j \Tr (\ad_{E_j}^* \ad_{E_j} S K) = - \sum_{ij} \<[E_i,E_j], [S K E_i, E_j]\>$.

Now for $\psi \in \Der(\l)$, let $\psi=S+K$, with $S$ symmetric and $K$ skew-symmetric. Then
$[\psi, \psi^*] = KS + (KS)^*$, so by \eqref{eq:sk},
$\Xi([\psi, \psi^*], \, \id)= 2 \, \Xi(KS, \,\id) = 2 \, \Xi(S,K) = 2 \, \Xi(S,\psi - S) = -2 \, \Xi(S,S)$.
This, together with \eqref{eq:riccidefdelta}, implies \eqref{eq:him}. Equation \eqref{eq:hinil} then follows from
\eqref{eq:tracestandard} and the fact that $B = 0$ and $H = 0$.

3. Let $\d_{ii}, \; i=1,2,$ be the set of all $\psi \in \End(\l)$ such that
$\psi(\l_i) \subset \l_i,\; \psi_{|\l_i} \in \Der (\l_i)$, and $\psi_{|\l_j} = 0$, for $j \ne i$. Clearly,
$\d_{ii} \subset \Der (\n)$. For $i \ne j$, let $\d_{ij}$ be the set of all $\psi \in \End(\l)$ such that
$\psi_{|\l_j \oplus [\l_i, \l_i]} = 0$ and $\psi(\l_i)$ lies in the center of $\l_j$. Then
$\Der (\l) = \oplus_{i,j=1}^2 \d_{ij}$.

Let now $\phi \in \d_{11} \oplus \d_{22}$ be the derivation of $\n$ such that $\phi_{|\l_i} = \phi_i$. Clearly,
$\phi$ is semisimple and real. Moreover, the equation \eqref{eq:pEtrace} holds for any
$\psi \in \d_{11} \oplus \d_{22}$, as each of the $\phi_i$'s is pre-Einstein, and for any
$\psi \in \d_{12} \oplus \d_{21}$, as both sides vanish.

The claim now follows from assertion (b) of Proposition~\ref{p:pEexistunique}.
\end{proof}

As to assertion 3, note that if the $\l_i$'s are Einstein nilradicals, with $\Phi_i$ the Einstein derivations,
then $\l$ is also an Einstein nilradical whose Einstein derivation is $\Phi_1 \oplus \Phi_2$ \cite[Theorem 4]{P}
(geometrically
this says that the Riemannian product of two standard Einstein solvmanifolds of the same Ricci curvature
is again a standard Einstein solvmanifold).

\begin{proof}[Proof of Theorem~\ref{th:pE}]
Let $\g$ be a solvable nonunimodular  Lie algebra, with an Einstein inner product $\<\cdot , \cdot \>$. Let $\n$ be the
nilradical of $\g$, and let $\ag = \n^{\perp}$. Denote $\lambda_1 < \ldots < \lambda_m$ the eigenvalues of
the pre-Einstein derivation $\phi = \phi_\n$ of $\n$, and $\n_1, \ldots, \n_m$ the corresponding eigenspaces.

Suppose that $\phi$ satisfies \eqref{eq:pE>0}, that is, $\lambda_i \ge 0$ and for any
$\psi \in \Der(\n), \quad \psi(\n_i) \subset \oplus_{j=i}^m \n_j$ (which is equivalent to $\ad_\phi \ge 0$).

Define $F \in \End(\g)$ by $F_{|\n} = \phi^{\perp}$ (as in assertion 1 of Lemma~\ref{l:gs}) and $F_{|\ag} = 0$.
Then $F$ is symmetric and nonnegative.
The eigenspaces of $F_{|\n}$ are $\m_i=\oplus_{j=i}^m \n_j \cap (\oplus_{j=i+1}^m  \n_j)^{\perp}$
(as constructed in Lemma~\ref{l:gs}), with the corresponding eigenvalues $\la_i$ in the increasing order.
Moreover, for any derivation $\psi$ of
$\n, \quad \psi(\m_i) \subset \psi(\oplus_{j=i}^m \n_j) \subset \oplus_{j=i}^m \n_j = \oplus_{j=i}^m \m_j$.
It follows that $\Tr (\psi \circ \phi^{\perp}) = \Tr (\psi \circ \phi) = \Tr \psi$, by definition of $\phi$.
As $F$ is symmetric, we obtain that for any $\Psi \in \Der(\g)$,
\begin{equation}\label{eq:trPsi}
\Tr \, (S(\Psi) \circ F) = \Tr \, (\Psi \circ F) = \Tr \, (\Psi_{|\n} \circ \phi^\perp)
= \Tr \, \Psi_{|\n}= \Tr \, \Psi= \Tr \, S(\Psi).
\end{equation}
Also, $\Tr \, F = \Tr \, \phi$, and $\Tr (F^2) = \Tr (\phi^2) = \Tr \, \phi$ (as $F$ and $\phi$ have the same
nonzero eigenvalues).

Substituting $F$ as $A$ to \eqref{eq:riccidefein} we find that the right-hand side is nonnegative. To see that,
choose an orthonormal basis $\{E_i\}$ to be of the form $\{Y_1, \ldots, Y_a, E_1, \ldots, E_n\}$, where
$E_s \in \n$ are the eigenvectors of $F_{|\n}$, and $Y_k \in \ag$.

Note that for any $E \in \m_i$ and $Y \in \ag, \quad \ad_YE \in \oplus_{j=i}^m \m_j$, so
$\<F[E, Y] - [FE, Y] - [E, FY], [E, Y]\> = \<[F, \ad_{Y}] E, \ad_{Y} E\> =
\<(F - \la_i \id)_{|\oplus_{j=i}^m \m_j} \ad_{Y} E, \ad_{Y} E\> \ge 0$.
This implies that
$\sum_{s,k} \<F[E_s, Y_k] - [FE_s, Y_k] - [E_s, FY_k], [E_s, Y_k]\> \ge 0$.

Furthermore,
$\sum_{s,t} \<F[E_s, E_t] - [FE_s, E_t] - [E_s, FE_t], [E_s, E_t]\> \ge 0$ by Lemma~\ref{l:gs}.
The remaining sum on the right-hand side is
$\sum_{k,l} \<F[Y_k, Y_l], [Y_k, Y_l]\> \ge 0$ as $F \ge 0$.

On the left-hand side of \eqref{eq:riccidefein}, the third summand vanishes as $F(\g) \subset \n$, so we get
$|c| \Tr F - \Tr (S(\ad_H) F) \le 0$. Substituting $A = \ad_H$ into \eqref{eq:riccidefein} we obtain
$|c| = \Tr \,S(\ad_H)^2/ \Tr \, S(\ad_H)$, so
$$
\Tr (S(\ad_H)^2) \, \Tr F \le \Tr (S(\ad_H) \circ F) \, \Tr S(\ad_H).
$$
As $\Tr F = \Tr (F^2) (=\Tr \phi)$ and $\Tr (S(\ad_H) \circ F) = \Tr \, S(\ad_H)$ from \eqref{eq:trPsi}, the
latter inequality gives
$\Tr (S(\ad_H)^2) \, \Tr (F^2) \le (\Tr (S(\ad_H) \circ F))^2$. By the Cauchy-Shwartz inequality, all the
inequalities above must be equalities. In particular, $\<F[Y_k, Y_l], [Y_k, Y_l]\> =0$, for all $Y_k, Y_l \in \ag$,
which implies that $\ag$ is abelian, as $F_{|\n} = \phi^\perp>0$.

\end{proof}

\begin{remark} \label{r:ad_phi}
As it can be seen from the proof, we need the condition $\ad_\phi \ge 0$ to be satisfied only on the subspace of
$\Der(\n)$ spanned by $\ad_{Y|\n}, \, Y \perp \n$. More precisely, it is sufficient to require that for each
$Y \perp \n, \; \ad_{Y|\n}$ lies in the direct sum of the eigenspaces of $\ad_\phi$ with nonnegative eigenvalues.
\end{remark}

\begin{proof}[Proof of Theorem~\ref{th:split}]
Let $(\g, \<\cdot,\cdot\>)$ be an Einstein metric solvable Lie algebra with the nilradical $\n$. Suppose that $\g$ is
split. Let $\h \subset \g$ be a Cartan subalgebra. Then $\g = \h + \n$ (not a direct sum, in general), as
$[\g,\g] \subset \n$. For any $X \in \h, \; (\ad_X)^S \h = 0$. As $\g$ is split, $(\ad_X)^S = \ad_Y$ for some
$Y \in \g$. Since $[\h, Y] = 0, \quad Y \in \h$, so $\h$ is also split. All the nilpotent elements from $\h$ lie in $\n$
(as $\n$ is the maximal nilpotent ideal of $\g$). It follows that $\h$ contains an abelian subalgebra $\ag$ such that
$\g= \ag \oplus \n$ and $\ad_Y$ is semisimple for any $Y \in \ag$. Let now $\phi$ be a pre-Einstein derivation for $\n$
and suppose that $\phi > 0$. By conjugation, we can choose $\phi$ to lie in the maximal reducible subalgebra of
$\Der (\n)$ containing the torus $\{\ad_{Y|\n}\, |\, Y \in \ag\}$. As $\phi$ lies in the center of that subalgebra
(see the proof of Proposition~\ref{p:pEexistunique}), it
commutes with all the $\ad_{Y|\n}, \, Y \in \ag$. As $\g = \ag \oplus \n$, any vector
$Z \perp \n$ differs from some $Y \in \ag$ by a vector $X \in \n$, so $\ad_{Z|\n} = \ad_{Y|\n} + \ad_{X|\n}$. Now
$\ad_{Y|\n} \in \Ker \ad_\phi$, as $\ad_{Y|\n}$ and $\phi$ commute. If $X = \sum_i X_i$, with $\phi X_i = \lambda_i X_i$,
then $\ad_\phi (\ad_{X_i|\n}) = \ad_{\phi(X_i)} = \lambda_i \ad_{X_i|\n}$. Since $\phi > 0$, it follows that
$\ad_{X|\n}$ (and therefore $\ad_{Z|\n}$) lies in the direct sum of the eigenspaces of $\ad_\phi$ with the nonnegative
eigenvalues. By Remark~\ref{r:ad_phi}, the proof can now be finished by the same arguments as in the proof of
Theorem~\ref{th:pE}.
\end{proof}

We end this section with two examples showing the limitations of Theorem~\ref{th:pE}.
Of course, if a nilpotent algebra $\n$ has no positive derivations, Theorem~\ref{th:pE} says nothing.
Example~\ref{ex:nonpositive} shows that the things could be worse: even when a positive derivation exists, the
pre-Einstein derivation can be nonpositive. Example~\ref{ex:adphinonpositive} shows that the first condition of
Theorem~\ref{th:pE} ($\phi >0$) does not imply the second one ($\ad_\phi \ge 0$).

\begin{example} \label{ex:nonpositive}
Any two-step nilpotent Lie algebra admits a positive derivation, the canonical derivation $\Phi$
($\Phi$ acts as $2 \, \id$ on the derived algebra and as $\id$ on its linear complement; see Section~\ref{s:2step} for
details). However, the pre-Einstein derivation of a two-step nilpotent Lie algebra can be nonpositive.

Let $\b = \oplus_{i=-1}^1 \b_i,\; \m = \oplus_{k=-2}^2 \m_k, \; \n = \b \oplus \m$ be the linear space decompositions.
Denote
$q_i = \dim \b_i$, $p_k = \dim \m_k$. Define the structure of a two-step nilpotent Lie algebra on $\n$ by requiring that
$[\b_i, \b_j] = \m_{i+j}$ and $[\n, \m] = 0$. Moreover, suppose that every two-step
nilpotent subalgebra $\b_i \oplus \m_{2i} \subset \n$, $i =-1,0,1$, admits only one semisimple derivation, up to
conjugation and scaling: the canonical derivation
(by \cite[Proposition 3.4.3]{E2}, for $p, q$ large enough and satisfying $2 < p < q(q-1)/2-2$, a generic two-step
nilpotent Lie algebra of type $(p,q)$ has that property).

Let $\psi_0$ be the derivation of $\n$ acting by multiplying every vector from $\b_i$ by $i$, and every vector
from $\m_k$ by $k$. For $a, b \in \mathbb{R}$, define $\psi_{a,b} = a \Phi + b \psi_0 \in \Der(\n)$.
We claim that one of the derivations $\psi_{a,b}$ is a pre-Einstein derivation for $\n$. Clearly,
every $\psi_{a,b}$ is semisimple and real. To check equation \eqref{eq:pEtrace}, it is sufficient to take only those
$\psi \in \Der(\n)$ which commute with $\phi$ and are semisimple. The restriction of any such $\psi$ to every subalgebra
$\b_i \oplus \m_{2i}$ must be a semisimple derivation, hence must be proportional to the canonical derivation of
$\b_i \oplus \m_{2i}$. It follows that $\psi = \psi_{a,b}$ for some $a, b \in \mathbb{R}$ (in fact,
$T=\{\psi_{ab} \, | \, a, b \in \mathbb{R}\}$ is a maximal torus of derivations of $\n$). Now, a derivation
$\psi_{x,y} \in T$ is pre-Einstein, if $\Tr \psi_{x,y} \psi_{a,b} = \Tr \psi_{a,b}$, for all $a, b \in \mathbb{R}$.
Solving this, we find that one of the eigenvalues of $\psi_{x,y}$ is
$x-y = (Q_1+(8 p_{-2} + 3 p_{-1}+ 4 q_{-1}+q_{0}) q_{1}
+(8 p_{-1}+2 p_{0}+2 p_2+9 q_{-1}+2 q_{0}+18 p_{-2}) p_{1}-p_{1} q_{1})/Q_2$,
where $Q_1$ and $Q_2$ are quadratic forms in the $q_i$ and $p_k$ with nonnegative coefficients, and $Q_1$ does not
contain $p_1$ and $q_1$. This expression can be made negative by choosing $p_1$ and $q_1$ large enough (the inequalities
$p_{i+j} \le q_iq_j$ and $2 < p_{2i} < \frac12 q_i(q_i-1)-2$ needed to have $[b_i, b_j] = \m_{i+j}$ and for every
$\b_i\oplus \m_{2i}$ to be generic can be easily achieved by multiplying all the $q_i$ and $p_k$ by a large natural
number).


\end{example}

\begin{example} \label{ex:adphinonpositive}
This example shows that the positivity of the pre-Einstein derivation $\phi$ of a nilpotent Lie algebra does not
imply the nonnegativity of $\ad_\phi$. Let $\n$ be a two-step nilpotent Lie algebra attached to the graph
$\mathcal{G}_{n,n,0}, \; n > 2$ (see \cite[Section 5]{LW}). The graph $\mathcal{G}_{n,n,0}$ has $2n+2$ vertices
$v_i$, its edge set is $E=\{e_{12}, e_{13}, \ldots, e_{1,n+2}, e_{2,n+3}, \ldots, e_{2,2n+2}\}$, where $e_{ij}$ is
the edge joining $v_i$ and $v_j$. The algebra $\n$, as a linear space, is spanned by the $v_i$'s and the
$e_{ij} \in E$, with the commutator relations given by $\z(\n) = \Span(E)$, and for $i < j, \; [v_i, v_j] = e_{ij}$, if
$e_{ij} \in E$, and $0$ otherwise. A direct computation shows that the pre-Einstein derivation $\phi$ of $\n$ has
the form $\phi \, v_i = x_i v_i,\; \phi \, e_{ij} = (x_i+x_j) \, e_{ij}$, up to conjugation, where
$x_1 = x_2 = 4 (n+2) d, \; x_i = (n+2) (n+4) d, \; i \ge 3$, and $d = ((n+4)^2-4)^{-1}$. It follows that $\phi > 0$
(one can show that this is always true for the algebras attached to a graph). For $\psi \in \Der(\n)$, which sends
$v_3$ to $e_{12}$ and all the other $v_i$'s and $e_{ij}$'s to zero, $\ad_\phi \psi = (2x_1 - x_3) \psi$, with the
coefficient $2x_1 - x_3 = (8+2n-n^2)\, d$ being negative for $n>4$.
\end{example}

\section{Applications of Theorem~\ref{th:pE}}
\label{s:appl}

In this section, we prove Theorem~\ref{th:nil} and Theorem~\ref{th:free}
applying Theorem~\ref{th:pE} to Einstein nilradicals and to free nilpotent Lie algebras, respectively.

\begin{proof}[Proof of Theorem~\ref{th:nil}]
To deduce Theorem~\ref{th:nil} from Theorem~\ref{th:pE}, it suffices to show that for an Einstein nilradical
the pre-Einstein derivation $\phi$ is (positively proportional to) the Einstein derivation and then to check that the
conditions of Theorem~\ref{th:pE} hold.

Let $\n$ be an Einstein nilradical, with $\Phi$ the Einstein derivation. Then $\Phi$ is semisimple, real and
satisfies \eqref{eq:tracestandard}, so $\phi = (-c)^{-1} \Phi$ is a pre-Einstein derivation by Definition~\ref{d:pE}.
As $\Phi$ is positive, $\phi$ is also positive. It remains to check that $\ad_\phi \ge 0$. Assume that there exists
$\psi \in \Der(\n)$ belonging to an eigenspace of $\ad_\phi$ with a negative eigenvalue. As $\phi = (-c)^{-1} \Phi$,
we get $[\Phi, \psi] = \lambda \psi,\; \lambda < 0$. Substituting such a $\psi$ into \eqref{eq:hinil} we obtain
$\lambda \Tr(\psi \psi^*) \ge 0$, which implies $\psi = 0$.
\end{proof}

Before giving the proof of Theorem~\ref{th:free}, we recall some facts about the structure of free Lie algebras and
their derivation algebras \cite{B}.
From among $p$-step nilpotent Lie algebras on $m \ge 2$ generators, a \emph{free Lie algebra} $\f(m, p)$ is the one
having the maximal dimension ($\f(m, p)$ is unique, up to an isomorphism). Given the generators
$e_1, \ldots, e_m$, the algebra $\f(m, p)$ is the linear span of all the $k$-folded brackets of the $e_i$'s, $k \le p$,
with the only relations between these brackets coming from the skew-symmetricity and the Jacobi identity. For every
$k=1, \ldots, p$, the subspace of $\f(m, p)$ spanned by the $k$-folded brackets is the space of \emph{Lie polynomials}
$\p(m, k)$. In particular, $\p(m, 1) = \Span (e_1, \ldots, e_m)$. The direct sum decomposition
$\f(m, p) = \oplus_{k=1}^p \p(m, k)$ is an $\mathbb{N}$-gradation. It corresponds to the \emph{canonical derivation}
$\Phi$ acting as a multiplication by $k$ on every $\p(m, k)$.

For a free Lie algebra, any assignment of the images to the generators extends to a (unique) derivation: for any linear
map $L: \p(m,1) \to \f(m,p)$, there exists a unique derivation whose restriction to $\p(m,1)$ coincides with $L$. In
particular, any endomorphism $L$ of $\p(m,1)$ extends to a derivation $\rho(L) \in \Der (\f(m,p))$ (for example,
$\rho(\id) = \Phi$).

\begin{proof}[Proof of Theorem~\ref{th:free}]
The spaces $\p(m,k)$ are
invariant with respect to $\rho(L)$. For every $k = 1, \ldots , p$, let $\rho_k(L)$ be the restriction of $\rho(L)$ to
the $\p(m,k)$. Then $\rho_k$ is a representation of the Lie algebra $\mathfrak{gl}(m)$ on the space $\p(m,k)$ of Lie
polynomials.

To deduce Theorem~\ref{th:free} from Theorem~\ref{th:pE}, it suffices to show that the pre-Einstein derivation
$\phi = \phi_{\f(m,p)}$ of the free Lie algebra $\f(m,p)$ is a positive multiple of $\Phi$. Then the condition
$\phi > 0$ is clearly satisfied. The condition $\ad_{\phi} \ge 0$ easily follows from the above, as for any
$\psi \in \Der (\f(m,p))$ and any $s =1, \ldots, p, \quad \psi(\p(m, s)) \subset \oplus_{k=s}^p \p(m, k)$.

As the pre-Einstein derivation is unique, up to conjugation (assertion (b) of Proposition~\ref{p:pEexistunique}),
to prove that
$\phi_{\f(m,p)}$ is indeed $\hat c \Phi$ (where $\hat c = \Tr \Phi \, (\Tr \Phi^2)^{-1} > 0$), it suffices to show
that $\hat c \, \Phi$ satisfies \eqref{eq:pEtrace}. The proof is literally the same as that of Lemma~8 of \cite{Ni}.
As $\Phi$ is semisimple, the derivation $\ad_\Phi$ of $\Der(\f(m,p))$ is also semisimple.
If $\psi \in \Der(\n)$ is an eigenvector of $\ad_\Phi$ with a nonzero eigenvalue, then
$\Tr \, \Phi \circ \psi = \Tr \, \psi = 0$, and \eqref{eq:pEtrace} is obviously satisfied. So it is
sufficient to consider only those $\psi$ which commute with $\Phi$. For any such $\psi$, the spaces $\p(m, k)$ are
invariant, and moreover, $\psi = \rho(L)$ for some endomorphism $L$ of $\p(m,1)$. As $\rho(\id) = \Phi$, it
suffices to prove the following: for any $L$ with $\Tr L = 0$ and for any $k \le p, \quad \Tr \rho_k(L) = 0$. This
follows from the fact that $\rho_{k|\mathfrak{sl}(m)}$ is a representation of the simple algebra $\mathfrak{sl}(m)$.
\end{proof}

\section{Nilradical having less than four generators}
\label{s:fili}

In this section, we establish Theorem~\ref{th:filiform} by proving a slightly more general fact: if
$\n$ is a nilpotent Lie algebra generated by less than four elements, then any Einstein metric solvable Lie algebra
$(\g, \<\cdot, \cdot \>)$ with $\g \supset \n \supset [\g,\g]$ is standard.

\begin{proof}
Let $\n$ be a nilpotent Lie algebra, and $\mathcal{C}_0 = \n,\; \mathcal{C}_p = [\mathcal{C}_{p-1}, \n], \; (p \ge 1)$
be the descending central series for $\n$. Denote $\n^1 := \n/\mathcal{C}_1$, with $\pi:\n \to \n^1$ the natural
projection.

Every endomorphism $A$ of $\n$ mapping $\mathcal{C}_1$ to itself (in particular, any derivation of $\n$) induces a
well-defined endomorphism $\pi (A): \n^1 \to \n^1$. In fact, $\pi$ is a homomorphism of the associative algebras, so
for any
$\psi_1, \psi_2 \in \Der(\n),\quad \pi(\psi_1\psi_2) = \pi(\psi_1) \pi(\psi_2)$. In particular, if $\psi_1$ commutes with
$\psi_2$, then $\pi(\psi_1)$ and $\pi(\psi_2)$ also commute, and if $\psi$ is a nilpotent
derivation, then $\pi(\psi)$ is nilpotent. The converse is also true: if $\pi(\psi)^m = 0$, then
$\psi^m (\n) \subset [\n, \n]$, and so $\psi^{mp+m-1} (\mathcal{C}_p) \subset \mathcal{C}_{p+1}$. Furthermore, if for a
derivation $\psi$ all the eigenvalues of $\pi(\psi)$ are real, then the same is true for $\psi$ itself. Indeed,
$\psi^S$, the semisimple part of $\psi$, is again a derivation, with the same eigenvalues as $\psi$. As $\psi^S$ is
semisimple, the endomorphism $\pi(\psi^S)$ is also semisimple. Moreover, $\psi-\psi^S$ is nilpotent and commutes with
$\psi$, which implies that $\pi(\psi)-\pi(\psi^S)$ is nilpotent and commutes with $\pi(\psi)$. It follows that
$\pi(\psi^S)$ is the semisimple part of $\pi(\psi)$. In particular, the eigenvalues of $\pi(\psi^S)$ are the same as
those of $\pi(\psi)$. Assume all of them are real. Choose the vectors
$X_1, \ldots, X_d \in \n \setminus \mathcal{C}_1$ such that $\pi (X_i)$ is the basis of eigenvectors for $\pi(\psi^S)$.
Every term $\mathcal{C}_p$ of the descending central series is spanned by $\mathcal{C}_{p+1}$ and all the $(p+1)$-folded
brackets $X_{i_1i_2\ldots i_{p+1}} = [\ldots[[X_{i_1}, X_{i_2}], X_{i_3}], \ldots ],X_{i_{p+1}}]$, where
$i_j \in \{1, \ldots, d\}$. Choosing, for every $p \ge 1$, a linearly independent set of the
$X_{i_1i_2\ldots i_{p+1}}$'s whose span complements $\mathcal{C}_{p+1}$ in $\mathcal{C}_p$ we get a basis for $\n$
having the property that every $X_{i_1i_2\ldots i_{p+1}}$ from that basis is an eigenvector of $\psi^S$ modulo
$\mathcal{C}_{p+1}$, with a real eigenvalue. It follows that all the eigenvalues of $\psi^S$ are real (the matrix of
$\psi^S$ is triangular, with a real diagonal).

Now let $(\g, \<\cdot, \cdot \>)$ be an Einstein metric solvable Lie algebra, $\n$ be a nilpotent ideal of $\g$
containing the derived algebra $[\g, \g]$, and let $\ag$ be the orthogonal complement to $\n$ in $\g$.

Define a linear map $\T: \ag \to \End(\n^1)$ by $\T(Y) = \pi(\ad_{Y|\n})$. Then $\T(\ag)$ is a commuting family
of endomorphisms of $\n^1$. For every $Y \in \Ker \T,\; \ad_Y$ is a nilpotent derivation of $\g$, hence
$Y \in \Ker B$, the kernel of the Killing form of $\g$. If $\dim \mathrm{Range}(\T) =1$, there is nothing to prove,
as then, by the Lie's Theorem, the nilradical of $\g$ has codimension $1$ and the claim follows from
Lemma~\ref{l:nilradvsderived}.

Assume now that $\dim \n^1 \le 3$. Clearly $\dim \n^1 \ne 1$. We will use \cite[Lemma 4.7]{H} saying that
any Einstein inner product on $\g$ is standard, provided the index of the Killing form $B$ of $\g$ is at most one.

Let $\dim \n^1 =2$. As $\mathrm{Range}(\T)$ is a commuting family of endomorphisms, $\dim \mathrm{Range}(\T) \le 2$,
with the equality only when there exists $Y \in \ag$ such that $\T(Y) = \id_{\n^1}$. It follows that the Killing form
$B$ of $\g$ has rank at most $2$ and there exists a vector $Y \in \ag$ such that $B(Y, Y) > 0$. Then
the index of $B$ is at most one, hence $(\g, \<\cdot, \cdot \>)$ is standard.

Let $\dim \n^1 = 3$. If all the eigenvalues of all the $\T(Y)$'s are real, then all the $\ad_Y$'s have only real
eigenvalues, so  $\g$ is completely solvable, hence standard (the index of $B$ is zero). If
one of the $\T(Y)$ has a nonreal complex eigenvalue, then $\T(Y)$ is semisimple for all $Y \in \ag$, and there is a
subspace $\ag'$ of codimension $1$ in $\ag$ such that all the eigenvalues of $\T(Y)$ are real, when $Y \in \ag'$.
Then the index of the Killing form $B$ of $\g$ is at most one, and again, any Einstein inner product on $\g$ is
standard.
\end{proof}

\section{Two-step nilradical}
\label{s:2step}

In this section, we consider several classes of two-step nilpotent Lie algebras and show that neither of them
contains the nilradical of a nonstandard Einstein metric solvable Lie algebra.
We start with some preliminary facts, mostly following \cite{E2}.

\subsection{Preliminaries}
\label{ss:2steppre}

A two-step nilpotent Lie algebra $\n$ of dimension $p + q$, is said to be \emph{of type} $(p,q)$, if its derived
algebra $\m = [\n, \n]$ has dimension $p$. Clearly, $\m \subset \z(\n)$, the center of $\n$, and
$1 \le p \le D:=\frac12 q (q-1)$.

Choose a subspace $\b$ complementary to $\m$ in $\n$ (in the presence of an inner product, we usually take
$\b = \m^\perp$). The Lie bracket defines (and is defined by) a skew-symmetric bilinear map
$J: \b \times \b \to \m$. For any $f \in \m^*$, the two-form $J_f \in \Lambda^2\b$ is defined by
$J_f(X, Y) = f([X, Y])$, for $X, Y \in \b$. If a particular basis $\{Z_k\}$ for $\m$ is chosen, we abbreviate
$J_{Z_k^*}$ to $J_k$. In the presence of an inner product, we identify $\m$ with $\m^*$ and write $J_Z$ for the
skew-symmetric operator on $\b$ defined by $\<J_ZX, Y\> = \<Z, [X, Y]\>, \; X, Y \in \b,\, Z \in \m$.

Choosing the bases $\{X_i\}$ for $\m$ and $\{Z_k\}$ for $\b$, we get a family of
skew-symmetric $q \times q$ matrices $J_1, \ldots, J_p$, which are linearly independent, as $\m = [\n, \n]$.
With respect to the choice of the bases, the matrices $J_i$ are defined up to a simultaneous transformation
$J_k \to GJ_kG^t, \; G \in \mathrm{SL}(q)$ and linear combinations $J_k \to \sum_r c_{kr} J_r$ with a nonsingular
matrix $(c_{kr})$. Every choice of the basis $\{X_i\}$ defines a point $W = \Span(J_1, \ldots , J_p)$ on the
Grassmannian
$G(p, \og(q))$ of $p$-planes in $\og(q)$. The isomorphism classes of two-step nilpotent Lie algebras of type $(p, q)$ are
in one-to-one correspondence with the orbits $[W]$ of the action of $\mathrm{SL}(q)$ on $G(p, \og(q))$ defined by
$\Span(J_1, \ldots , J_p) \to \Span(GJ_1G^t, \ldots , GJ_pG^t), \; G \in \mathrm{SL}(q)$. The space $X(p, q)$ of such
orbits (the space of the isomorphism types of the two-step nilpotent Lie algebras of type $(p, q)$) is in general
non-Hausdorff.
The dimension count suggests that the stabilizer $\mathrm{SL}(q)_W$ of a generic $W \in G(p, \og(q))$ (``the stabilizer
in general position") is finite. This is indeed the case, with a few exceptions, as it is proved
in \cite{E1}. 

The splitting $\n = \b \oplus \m$ of a two-step nilpotent Lie algebra $\n$ is a gradation, which corresponds to the
\emph{canonical derivation} $\Phi$ defined by $\Phi (X + Z) = X + 2 Z$, for any $X \in \b,\; Z \in \m$.
The algebra of derivations $\Der(\n)$ of $\n = \b \oplus \m$ splits into a semidirect sum of the abelian ideal
$\mathfrak{J} = \{\psi \in \Der(\n)\,|\, \psi(\b) \subset \m, $ $\psi(\m) = 0\}$ consisting of
nilpotent derivations and the subalgebra $\mathfrak{G} = \mathbb{R}\Phi \oplus \mathfrak{sl}_W$, where
$\mathfrak{sl}_W$ is the Lie algebra of the stabilizer $\mathrm{SL}(q)_W$ of $W \in G(p, \og(q))$
\cite[Proposition 3.4.5]{E2}. Clearly, $\Tr\, \psi = 0$ for every $\psi \in \mathfrak{sl}_W$.

In the fixed bases for $\b$ and $\m$, any
$\psi \in \Der(\n)$ is represented by a matrix of the form
$\left(\begin{smallmatrix} F & 0 \\U & M \end{smallmatrix}\right)$, where $U$ is an arbitrary $p \times q$-matrix
(the set of the $\left(\begin{smallmatrix} 0 & 0 \\U & 0 \end{smallmatrix}\right)$'s is $\mathfrak{J}$),
and the $q \times q$-matrix $F$ and the $p \times p$-matrix $M$ satisfy
  \begin{equation}\label{eq:der2step}
    g(F)J_k:= J_kF + F^tJ_k = \sum\nolimits_{r=1}^p M_{kr} J_r.
  \end{equation}
For a given endomorphism $F$, the corresponding $M:= M(F)$ is defined by equation \eqref{eq:der2step}
uniquely, as the $J_r$'s are linearly independent.
With any $\psi =\left(\begin{smallmatrix} F & 0 \\U & M \end{smallmatrix}\right) \in \Der(\n)$,
each of the following is also a derivation:
$$
\left(\begin{matrix} 0 & 0 \\U & 0 \end{matrix}\right), \;
\left(\begin{matrix} F & 0 \\0 & M \end{matrix}\right), \;
\left(\begin{matrix} F^S & 0 \\0 & M^S \end{matrix}\right), \;
\left(\begin{matrix} F^N & 0 \\0 & M^N \end{matrix}\right), \;
\left(\begin{matrix} F^\mathbb{R} & 0 \\ 0 & M^{\mathbb{R}} \end{matrix}\right),
$$
where the superscripts $S, N$ and $\mathbb{R}$ denote the semisimple part, the nilpotent part, and the real
part, respectively. It is easy to see that $M(F^S) = M(F)^S, \; M(F^N) = M(F)^N$, and
$M(F^{\mathbb{R}}) = M(F)^{\mathbb{R}}$.

Let $\g$ be a solvable Lie algebra with the two-step nilradical $\n = \b \oplus \m$, and let $\ag$ be a linear space
complementary to $\n$. For every $Y \in \ag$, consider the $\mathfrak{G}$-part of $\ad_{Y|\n}$ and define the
endomorphisms $F(Y) = \pi_\b \circ \ad_{Y|\b}:\b \to \b$ and $M(Y) = M(F(Y)) = \ad_{Y|\m}:\m \to \m$ (where $\pi_\b$
is the linear projection to $\b$). Then both $\{F(Y)\, |\, Y \in \ag\}$ and $\{M(Y)\, |\, Y \in \ag\}$ are commuting
families of endomorphisms, and the same is true for their semisimple parts and real parts.

Let $\n = \b \oplus \m$ be a two-step nilpotent Lie algebra whose isomorphism type is defined by $[W]$, the
$\mathrm{SL}(q)$-orbit of the point $W \in G(p, \og(q))$. Let $W^\perp \in G(D-p, \og(q))$ be the $(D-p)$-dimensional
subspace of $\og(q)$ orthogonal to $W$ with respect to the inner product $\<J,K\> = - \Tr JK$ on $\og(q)$.
The two-step nilpotent Lie algebra $\overline{\n}=\b \oplus \overline{\m}, \; \dim \overline{\m} = D - p$
defined by $[W^\perp]$ is called \emph{dual} to $\n$.
It is easy to see that $\overline{\n}$ is well-defined and that $\overline{\overline{\n}}$ is isomorphic to $\n$.
Moreover, as $\mathrm{SL}(q)_W$ is isomorphic to $\mathrm{SL}(q)_{W^\perp}$ ($\mathrm{SL}(q)_{W^\perp}$ consists
of the matrices transposed to those from $\mathrm{SL}(q)_W$),
to every $\psi =\left(\begin{smallmatrix} F & 0 \\0 & M(F) \end{smallmatrix}\right) \in \Der(\n)$, there corresponds
the derivation $\overline{\psi} =
\left(\begin{smallmatrix} F^t & 0 \\0 & \overline{M}(F^t) \end{smallmatrix}\right) \in \Der(\overline{\n})$.

\subsection{Two-step nilradicals of class $\Og(p,q)$}
\label{ss:2stepopq}

In this section, we prove assertion 1 of Theorem~\ref{th:twostepopen} (Lemma~\ref{l:opqstandard} and
Lemma~\ref{l:opqlarge}).

Let $\n$ be a two-step nilpotent Lie algebra of type $(p,q)$ with $\m = [\n, \n]$.
Consider two linear forms, $t$ and $t_1$, on $\Der(\n)$ defined by $t(\psi):= \Tr \psi, \; t_1(\psi):= \Tr \psi_{|\m}$.

\begin{definition}\label{d:opq}
$\Og(p,q)$ is the set of two-step nilpotent Lie algebras of type $(p,q)$ for which $\Ker t \subset \Ker t_1$.
\end{definition}


We begin with the following easy observation.

\begin{lemma}\label{l:opqstandard}
1. A two-step nilpotent Lie algebra $\n$ belongs to $\Og(p,q)$ if and only if
its pre-Einstein derivation is a multiple of the canonical derivation $\Phi$:
$$
\phi_\n = \mu \Phi, \quad \text{where} \quad \mu = \Tr \Phi / \Tr \Phi^2 = (q+2p)/(q+4p).
$$

2. Any Einstein metric solvable Lie algebra with the nilradical $\n \in \Og(p,q)$ is standard.
\end{lemma}

\begin{proof}
1. This follows from Definition~\ref{d:pE} and from the uniqueness of the pre-Einstein derivation (assertion
(b) of Proposition~\ref{p:pEexistunique}).

2. This follows from Theorem~\ref{th:pE}: both conditions $\Phi > 0$ and $\ad_\Phi \ge 0$ are clearly satisfied.
\end{proof}

The next lemma shows that the set $\Og(p,q)$ is ``massive":

\begin{lemma}\label{l:opqlarge}
1. If $\n \in \Og(p,q)$, then $\overline{\n} \in \Og(D-p,q)$.
If $\n_i \in \Og(p_i,q_i)$ for $i =1, 2$, then the direct sum of the Lie algebras $\n_1$ and $\n_2$
lies in $\Og(p_1+p_2, q_1+q_2)$.

2. For all the pairs $(p, q)$ such that $1 < p < D-1$, except when $q$ is odd and $p=2$ or $p=D-2$, the set
$\{W \,|\, [W] \in \Og(p,q)\}$ contains an open and dense subset of $G(p, \og(q))$.
\end{lemma}

\begin{remark} Although the direct sum of a two-step Lie algebra $\n \in \Og(p,q)$ and an abelian ideal $\mathbb{R}^m$
does not belong to $\Og(p,q+m)$, it is not difficult to see that any Einstein metric solvable Lie algebra having
$\n \oplus \mathbb{R}^m$ as the nilradical ($\n \in \Og(p,q)$) is standard using Theorem~\ref{th:pE} and
assertion 3 of Lemma~\ref{l:gs}.
\end{remark}

\begin{remark} \label{r:missing}
The pairs $(p, q)$ not covered by assertion 2 of Lemma~\ref{l:opqlarge} are $(1, q),\; (D-1, q),\; (D, q)$, and
$(2, 2k+1), \; (D-2, 2k+1)$. There is only one two-step nilpotent Lie algebra of type $(D, q)$, the free algebra. It is
an Einstein nilradical by \cite[Proposition 2.9 (iii)]{GK}. If $p=1$, the algebra $\n$ is isomorphic to the
direct sum of a Heisenberg algebra and an abelian ideal, hence is an Einstein nilradical. Every algebra of type
$(D-1, q)$ is an Einstein nilradical, as it is proved in Lemma~\ref{l:p=D-1} below. By Theorem~\ref{th:nil}, neither of
the above algebras can be the nilradical of a nonstandard Einstein metric solvable Lie algebra.

In contrast, the two-step nilpotent Lie algebras of types $(2, 2k+1)$ and $(D-2, 2k+1)$ are ``the worst" from the
point of view of Definition~\ref{d:opq}: the sets $\Og(2,2k+1)$ and  $\Og(D-2,2k+1)$ are \emph{empty}. Indeed, from
\cite[Theorem 5.1]{LR} it follows that
any two odd-dimensional skew-symmetric matrices have two common complementary isotropic subspaces of nonequal
dimensions. Therefore, for every two-step nilpotent Lie algebra $\n = \b \oplus \m$ of type $(2, 2k+1)$,
there is a decomposition $\b = \b_1 \oplus \b_2, \; \dim \b_i = q_i$ such that $[\b_1, \b_1]= [\b_2, \b_2]= 0$ and
$q_1 \ne q_2$. It follows that an endomorphism $\psi$ of $\n$, acting as the identity on $\b_1$, minus the identity on
$\b_2$, and zero on $\m$, is a derivation. Then the trace of the derivation $\eta = (2k+5)\psi + (q_2-q_1)\Phi$ is
zero, but $t_1(\eta) = \Tr \eta_{|\m} = 2(q_2-q_1) \ne 0$. It follows that $\Og(2,2k+1) = \varnothing$, and hence
$\Og(D-2,2k+1) = \varnothing$ by assertion 1 of Lemma~\ref{l:opqlarge}.
\end{remark}

\begin{proof}
1. For any
$\psi =\left(\begin{smallmatrix} F & 0 \\0 & M(F) \end{smallmatrix}\right) \in \Der(\n)$,
there exists (a unique) $\overline{\psi} \in \Der(\overline{\n})$ of the form
$\overline{\psi} =\left(\begin{smallmatrix} F^t & 0 \\0 & \overline{M}(F^t) \end{smallmatrix}\right)$.
In particular, $W^\perp \subset \og(q)$ (with respect to the inner product $\<J, K\> = - \Tr JK$ on $\og(q)$)
is an invariant subspace of the linear operator $g(F^t)$ acting on $\og(q)$.
The subspace $W \subset \og(q)$ is in general not an invariant subspace of $g(F^t)$. Define
$\hat g(F^t): W \to W$ by $\hat g(F^t)= \pi_W \circ g(F^t)_{|W}$, where $\pi_W$ is the orthogonal projection to $W$.
Choosing an orthonormal basis $\{J_k\}$ for $W$ we easily obtain that
$\Tr \hat g(F^t) = -\sum_{k=1}^q \Tr ((J_kF^t + FJ_k)J_k) = -\sum_{k=1}^q \Tr ((J_kF + F^tJ_k)J_k) = \Tr g(F)_{|W}$.

As the trace of the linear operator $g(F^t)$ (acting on the whole $\og(q)$) equals
$\Tr g(F^t)_{|W^\perp} + \Tr \hat g(F^t)$, we obtain $\Tr g(F^t)_{|W^\perp} + \Tr g(F)_{|W} = \Tr g(F^t)$.
An easy computation shows that $\Tr g(F^t) = (q-1) \Tr F$. As
$\Tr g(F)_{|W} = \Tr M(F)$ and $\Tr g(F^t)_{|W^\perp}= \Tr \overline{M}(F^t)$, we find that
$\Tr M(F) + \Tr \overline{M}(F^t)= (q-1) \Tr F$, which implies $\Tr \psi + \Tr \overline{\psi} = (q+1) \Tr F$.

Assume now that $\Tr \overline{\psi} = 0$. Then $\Tr \psi = (q+1) \Tr F$, which implies $\Tr F = 0$, as
from $\n \in \Og(p,q)$ it follows that $\Tr \psi = (2p/q+1) \Tr F$ (and as $2p < q^2$).
This proves the first statement of assertion 1.

The second statement follows directly from assertion 3 of Lemma~\ref{l:gs} and assertion 1 of Lemma~\ref{l:opqstandard}.

2. According to \cite[Proposition 3.4.3]{E2}, for every pair $(p, q)$ such that $1 < p < D-1$, with few exceptions,
there exists an open and dense $\mathcal{O}(p,q) \subset G(p, \og(q))$ such that for any
$W \in \mathcal{O}(p,q)$, the two-step nilpotent algebra $\n$ corresponding to $W$ has the following property:
any automorphism of $\n$ with determinant $1$ has only eigenvalues whose module is $1$.

The excepted pairs are $(p,q) = (3, 4),\; (3, 5),\; (3, 6),\; (2, q)$ and the dual ones.

Clearly, if $W \in \mathcal{O}(p,q)$, then $[W] \in \Og(p,q)$, as for the corresponding two-step nilpotent Lie algebra
$\n$, any derivation whose trace vanishes has zero real part. This proves the assertion for all the ``generic" pairs.
Moreover, as for the dual pairs $(p, q)$ and $(D-p, q)$, both the Grassmannians $G(p,\og(q))$ and $G(D-p,\og(q))$ and
the spaces $\Og(p,q)$ and $\Og(D-p,q)$ are homeomorphic (assertion 1), it remains to consider only the following cases:
$(p, q) = (3, 4),\; (3, 5), \; (3, 6),\; (2, 2k)$.

Let $\mathcal{S} \subset (\og(q))^p$, for a given pair $(p, q)$, be the set of linearly independent $p$-tuples of
skew-symmetric $q \times q$ matrices. The projection $\Span: \mathcal{S} \to G(p,\og(q))$ sends open and dense subsets
to open and dense subsets. Note that to prove the assertion it suffices to show that the set
$\{W: [W] \in \Og(p,q)\} \subset G(p, \og(q))$ (or its $\Span^{-1}$ preimage in $\mathcal{S}$) has a nonempty interior.
Indeed,
by Definition \ref{d:opq}, the fact that $[W] \in \Og(p,q)$ for $W= \Span (J_1, \ldots, J_p)$ means that all the
solutions of a certain system of linear equations (expressing the fact that $\psi$ is a derivation whose trace is zero)
satisfies another linear equation, $t_1(\psi) = 0$. The entries of the matrix $\mathcal{M}$ of that system are linear
in the entries of the $J_k$'s. On an open and dense subset $\mathcal{S}' \subset \mathcal{S}$,
$\rk \, \mathcal{M}$ is maximal, and the condition that any solution satisfies $t_1(\psi) = 0$ is algebraic in the
entries of $\mathcal{M}$. If it is satisfied on some open subset of $\mathcal{S}'$, then it is satisfied on an open
and dense subset of $\mathcal{S}'$, and hence of $\mathcal{S}$.

\underline{Case $(p, q) = (2, 2n)$}. Let $\hat{\mathcal{S}} \subset \mathcal{S} \subset \og(2n)^2$ be the set of the
pairs $(J_1, J_2)$ such that $\det J_1J_2 \ne 0$, and the polynomial $\det (J_2 - x J_1)$ has no multiple roots for
$x \in \mathbb{C}$. Then $\hat{\mathcal{S}}$ is open and dense in $\mathcal{S}$ (being a nonempty complement to the zero
set of a certain system of polynomial equations). By \eqref{eq:der2step}, the fact that $[W] \in \Og(2,2n)$ for
$W = \Span(J_1, J_2)$ is equivalent to the following: any
two matrices $F$ and $M$ satisfying the equations $J_kF + F^tJ_k = \sum_{r=1}^2 M_{kr} J_r, \; k =1, 2$, also
satisfy the equation $\Tr \, F  = 2n \Tr \, M$.
This, in turn, is equivalent to the existence of a pair of matrices $K_1, K_2 \in \og(2n)$ such that
\begin{equation}\label{eq:belongtoopq}
K_1 J_1 + K_2 J_2 = 2 \, I_{2n}, \quad \Tr (K_rJ_s) = 2n \K_{rs}.
\end{equation}

For any pair $(J_1, J_2) \in \hat{\mathcal{S}}$ the system \eqref{eq:belongtoopq} has a solution. Indeed,
by \cite[Theorem 5.1]{LR}, we can simultaneously reduce the matrices $J_1, J_2$ to the form
$J_1 = \left(\begin{smallmatrix} 0 & I_n \\ -I_n & 0 \end{smallmatrix}\right), \;
J_2 = \left(\begin{smallmatrix} 0 & C \\ -C^t & 0 \end{smallmatrix}\right)$, where $C \in \mathrm{GL}(n)$. For
$K_r = \left(\begin{smallmatrix} 0 & T_r \\ -T_r^t & 0 \end{smallmatrix}\right)$, the system \eqref{eq:belongtoopq}
is equivalent to
$T_1 = -2I_n - T_2 C^t, \; [T_2, C^t] = 0,\; \Tr \, T_2 = 0$,
$\Tr \, (T_2 C^t) = -n, \; \Tr \, (T_2 C^{t2} + 2 C^t) = 0$. If $n > 2$, we can easily find a complex solution taking
$T_2$ diagonal in the Jordan basis for $C^t$ and using the fact that all the eigenvalues of $C^t$ are distinct. The
real part of it gives a real solution to the system \eqref{eq:belongtoopq}. If $n = 2$, take
$T_2 = (\la_1-\la_2)^{-2} ((2 \Tr \, C^t) I_2 - 4 C^t)$, $T_1 = -2I_2 - T_2 C^t$ (where $\la_1 \ne \la_2$
are the eigenvalues of $C^t$).

\underline{Case $(p, q) = (3, 4)$}. There exists a $W \in G(3, \og(4))$ such that the two-step nilpotent Lie algebra
$\n$ defined by $[W]$ is nonsingular, hence belongs to $\Og(3,4)$ by Lemma~\ref{l:nonsingular} below. The set of such
$W$'s has a nonempty interior, so $\{W \,|\, [W] \in \Og(3,4)\}$ contains an open and dense subset of
$G(3, \og(4))$ (we prove a stronger fact in assertion 1 of Lemma~\ref{l:small}).

\underline{Case $(p, q) = (3, 5)$}. As it follows from \cite[5.4]{E1}, the action of $\mathrm{SL}(5)$ on
$G(3, \og(5))$ has open orbits, and there exist locally rigid two-step nilpotent algebras of type $(3,5)$.
As $\dim \mathrm{SL}(5) = 24$ and $\dim G(3, \og(5)) = 21$, any $W \in G(3, \og(5))$ whose stabilizer is
three-dimensional defines a rigid algebra. A direct computation based on \eqref{eq:der2step} shows that
such a $W$ can be taken as the span of the matrices
\begin{equation*}
J_1 = \begin{pmatrix} 0 & 0 & 0 & 0 & 1 \\ 0 & 0 & 0 & 1 & 0 \\ 0 & 0 & 0 & 0 & 0
\\ 0 & -1 & 0 & 0 & 0 \\ -1 & 0 & 0 & 0 & 0\end{pmatrix} \quad
J_2 = \begin{pmatrix} 0 & 0 & 0 & 1 & 0 \\ 0 & 0 & 1 & 0 & 0 \\ 0 & -1 & 0 & 0 & 0
\\ -1 & 0 & 0 & 0 & 0 \\ 0 & 0 & 0 & 0 & 0\end{pmatrix} \quad
J_3 = \begin{pmatrix} 0 & 0 & 0 & 0 & 0 \\ 0 & 0 & 0 & 0 & 1 \\ 0 & 0 & 0 & 1 & 0
\\ 0 & 0 & -1 & 0 & 0 \\ 0 & -1 & 0 & 0 & 0\end{pmatrix}.
\end{equation*}
From \eqref{eq:der2step} we find that for any derivation of the two-step nilpotent Lie algebra $\n$ defined by $[W]$,
\begin{equation*}
F = x I_5 + \begin{pmatrix}
-2y & -2z & 0 & 0 & 0 \\
w & -y & -z & 0 & 0 \\
0 & 3w & 0 & 3z & 0 \\
0 & 0 & -w & y & z \\
0 & 0 & 0 & -2w & 2y\end{pmatrix} \qquad
M = 2x I_3 + \begin{pmatrix}
0 & -w & -z \\
z & -y & 0 \\
w & 0 & y \end{pmatrix},
\end{equation*}
so $\n \in \Og(3,5)$. As $\n$ is locally rigid, the set $\{W \,|\, [W] \in \Og(3,5)\}$ has a nonempty interior, and
hence contains an open and dense subset of the Grassmannian $G(3, \og(5))$.

\underline{Case $(p, q) = (3, 6)$}. We start with the following observation:
there exists an open subset $\mathcal{U} \subset G(3, \og(6))$ such that for any
$W = \Span(J_1, J_2, J_3) \in \mathcal{U}$, all the three $J_k$'s have a pair of complementary three-dimensional
isotropic
subspaces. To see that, consider a vector bundle $E \subset G(3, 6) \times \og(6)$ over the Grassmannian $G(3, 6)$ whose
fiber over a point $L \in G(3, 6)$ is the subspace $V_L \subset \og(6)$ of the ($6 \times 6$ skew-symmetric) matrices
having $L$ as an isotropic subspace. Let $G(3,E) \subset G(3, 6) \times G(3, \og(6))$ be the corresponding Grassmann
bundle: the fiber over $L \in G(3, 6)$ is the Grassmannian $G(3, V_L)$. Both $G(3,E)$ and $G(3, \og(6))$ are compact
differentiable manifolds of dimension $36$, and the projection on the second factor
$\pi: G(3,E) \to G(3, \og(6)), \; \pi (L, W) = W$, is a differentiable mapping. Take a point $(L, W) \in G(3,E)$ and
choose a basis $e_i$ for
$\mathbb{R}^6$ in such a way that $L = \Span (e_1, e_2, e_3)$ and $W = \Span (J_1, J_2, J_3)$, where
$J_k =\left(\begin{smallmatrix} 0 & T_k \\-T_k^t & C_k \end{smallmatrix}\right),\; k = 1, 2, 3$. It is not difficult to
see that $\Ker d\pi (L, W) = 0$, if there is no nonzero $3 \times 3$ matrix $S$ such that all the matrices
$T_kS$ are symmetric, and that all the $J_k$'s have a tree-dimensional isotropic subspace complementary to $L$, if
there exists a $3 \times 3$ matrix $P$ such that $T_k^tP - P^t T_k = C_k$, for all $k=1, 2, 3$. Both these conditions are
satisfied, provided the corresponding $9 \times 9$ determinants $d_1, d_2$ (whose entries are linear in the entries of
the $T_k$'s) do not vanish. A direct computation shows that $d_1 \ne 0$ and $d_2 \ne 0$, if we take say
$T_1 = I_3, \; T_2 = \diag(1,2,3)$, and $T_3$ skew-symmetric, with $(T_3)_{12}=1,\,(T_3)_{13}=2,\,(T_3)_{23}=3$.
It follows that the subset of the triples $(T_1, T_2, T_3)$, for which $d_1 \ne 0$ and $d_2 \ne 0$ is open and dense in
$\mathfrak{gl}(3)^3$. As the conditions $d_1, d_2 \ne 0$ do not depend of the choice of the basis for $W$, the subset
of those $W \in G(3, V_L)$, for which $\Ker d\pi (L, W) = 0$ and all the matrices from $W$ have a three-dimensional
isotropic subspace complementary to $L$, is open and dense in the fiber $G(3, V_L)$. On that subset, $\pi$ is locally
onto, so the Grassmannian $G(3, \og(6))$ contains an open subset $\mathcal{U}$ such that for any
$W = \Span(J_1, J_2, J_3) \in \mathcal{U}$, all the $J_k$'s have a pair of complementary three-dimensional isotropic
subspaces.

By \cite[5.4]{E1}, there exists an open and dense subset $\mathcal{O} \subset G(3, \og(6))$ such that the
stabilizer subgroup $\mathrm{SL}(6)_W$ of every $W \in \mathcal{O}$ has dimension one.
For any such $W$, the Lie algebra $\mathfrak{sl}(6)_W$ is one-dimensional, so for every two-step nilpotent Lie algebra
$\n$ of type $(3, 6)$ defined by some $[W], \; W \in \mathcal{O}$, we get
$\Der(\n) = \mathfrak{J} \oplus \mathbb{R}\Phi \oplus \mathbb{R}\psi_0$ for some derivation $\psi_0$ with
$\Tr \, \psi_0 = 0$.
Both linear forms $t$ and $t_1$ from Definition~\ref{d:opq} vanish on $\mathfrak{J}$, so
$\n \in \Og(3, 6)$ if and only if $t_1(\psi_0) = 0$. Now take
$W = \Span(J_1, J_2, J_3)$ from the open set $\mathcal{U} \cap \mathcal{O}$. For any $\n = \b \oplus \m$ defined by
such a $W$, we have $\b = \b_1 \oplus \b_2$, where $\b_1$ and $\b_2$ are three-dimensional isotropic subspaces of all
the three $J_k$'s. Then a semisimple endomorphism $\psi$ defined by
$\psi(X_1 + X_2 + Z) = X_1 - X_2, \; X_1 \in \b_1,  \; X_2 \in \b_2,  \; Z \in \m$, is a derivation of $\n$. As
$\Tr \, \psi = 0$, we can take $\psi_0 = \psi$. Since $\psi_{|\m} = 0,\; t_1(\psi) = 0$, so $\n \in \Og(3, 6)$. The
claim follows, as $\mathcal{U} \cap \mathcal{O}$ is a nonempty open subset of the Grassmannian $G(3, \og(6))$.
\end{proof}

\subsection{Two-step nonsingular nilradicals}
\label{nonsing}
Assertion 2(a) of Theorem~\ref{th:twostepopen}, which we prove in this section, follows from assertion 2 of
Lemma~\ref{l:opqstandard} and Lemma~\ref{l:nonsingular} below.

Let $\n$ be a two-step nilpotent Lie algebra with the derived algebra $\m$.
Algebra $\n$ is called \emph{nonsingular}, if for any $X \in \n \setminus \m$, the map
$\ad_n: \n \to \m$ is surjective \cite{E2}. The reason for the name comes from the following equivalent definition.
Choose an arbitrary inner product on $\n$ and set $\b = \m^\perp$. Then $\n$ is nonsingular
if and only if for any nonzero $Z \in \m$, the operator $J_Z$ is nonsingular.
The nonsingular two-step nilpotent Lie algebras are one of ``the most well behaved of all two-step nilpotent
Lie algebras" \cite{E2}. In particular, many interesting examples of Einstein nilradicals are two-step nonsingular,
for instance, the Heisenberg-type algebras \cite{GK}.

Note that the nonsingularity condition imposes strong restrictions on the type. The obvious one is that
$q=\dim \b$ must be even. A subtler dimension restriction comes from topology: for a (metric) nonsingular
two-step nilpotent Lie algebra,
the unit sphere $S^{q-1} \subset \b$ must admit $p$ continuous and pointwise linearly independent vector
fields defined (with respect to some basis $\{Z_k\}$ for $\m$) by $X \to J_{k}X$, for $X \in S^{q-1}$.
It follows from the Adams vector field Theorem, that $p \le \rho(q) - 1$, where $\rho(q)$
is the Radon-Hurwitz number.

\begin{lemma}\label{l:nonsingular}
Any nonsingular two-step nilpotent Lie algebra of type $(p,q)$ belongs to $\Og(p,q)$.
\end{lemma}
\begin{proof}
Let $\n$ be a nonsingular two-step nilpotent Lie algebra.
To show that $\n \in \Og(p,q)$, it suffices
to check that for any $\psi \in \Der(\n)$ of the form
$\psi =\left(\begin{smallmatrix} F & 0 \\0 & M \end{smallmatrix}\right)$, with both $F$ and $M$ semisimple and real,
the equation $\Tr \psi = 0$ implies $\Tr M = 0$ (or, equivalently, $\Tr F =0$).
Let $\la_k$ be the eigenvalues of $M$, with $Z_k$ the corresponding eigenvectors. According
to \eqref{eq:der2step} we have $J_k F + F^tJ_k= \la_k J_k$, for all $k = 1, \ldots, p$, so
$2J_k F = \la_k J_k + S_k$, where $S_k$ are some symmetric matrices.
It follows that $2F = \la_k I_q + J_k^{-1}S_k$, so $q \la_k = 2\Tr F$, for all $k$. This implies
$\Tr \psi = \Tr F + \Tr M = (1+2p/q)\Tr F$.
\end{proof}

\begin{example}
Note that for a two-step nilpotent Lie algebra $\n$ of type $(p, q)$, the fact that $\n \in \Og(p,q)$ and even a much
stronger condition of nonsingularity does not guarantee that $\n$ is an Einstein nilradical, as the following example
shows. Let $\n$ be a nonsingular two-step nilpotent Lie algebra of type $(2, 2n)$.
As it follows from Lemma~\ref{l:nonsingular} and assertion 1 of Lemma~\ref{l:opqstandard}, if $\n$ is an Einstein
nilradical, its eigenvalue type must be $(1< 2; 2n, 2)$,
so $\n$ is defined by an orbit $[W]$ of $W =\Span (J_1, J_2) \in G(2, \og(2n))$ such that $J_1^2+ J_2^2= -I_{2n}$
and $\det(xJ_1+yJ_2) \ne 0$ unless $x=y=0$. It follows that the matrix $K = J_2^{-1}J_1$ is normal, hence semisimple.
So for any $a, b, c, d \in \mathbb{R}$ with $ad-bc \ne 0$, the matrix
$(aJ_1+bJ_2)^{-1}(cJ_1+dJ_2) = (aK+bI_{2n})^{-1}(cK + bI_{2n})$ is also semisimple. This property is preserved under the
action of $\mathrm{GL}(2n)$ on $G(2, \og(2n))$ (which sends $J$ to $TJT^t$ for $T \in \mathrm{GL}(2n)$). Take now two
skew-symmetric $8 \times 8$ matrices of the form
$\tilde J_i =\left(\begin{smallmatrix} 0 & B_i \\-B_i^t & 0 \end{smallmatrix}\right)$, where $B_2 = I_4$,
$B_1 =\left(\begin{smallmatrix} J & I_2 \\0 & J \end{smallmatrix}\right)$, and
$J =\left(\begin{smallmatrix} 0 & 1 \\-1 & 0 \end{smallmatrix}\right)$. Clearly, any nontrivial linear
combination of the $\tilde J_i$'s is nonsingular, but $(\tilde J_2)^{-1}\tilde J_1$ is not semisimple.
\end{example}

\subsection{Two-step nilradicals of types $\mathbf{(D-1,q), (2, 5)}$ and $\mathbf{(3,4)}$}
\label{ss:D-1,34}

In this section, we consider two classes of two-step nilradicals whose type is not covered by assertion 1 of
Lemma~\ref{l:opqlarge} and prove assertions 2(b) and 2(c) of Theorem~\ref{th:twostepopen}.

\begin{lemma}\label{l:p=D-1}
Any two-step nilpotent Lie algebra of type $(D-1,q)$ is an Einstein nilradical.
\end{lemma}

As a consequence of Theorem~\ref{th:nil}, such an algebra cannot be the nilradical of any nonstandard Einstein metric
solvable Lie algebra.

\begin{proof}
The isomorphism type of a two-step nilpotent Lie algebra $\n$ of type $(D-1,q)$ is completely determined by the
one-dimensional subspace in $\og(q)$ orthogonal to $W$. What is more, if
$W^\perp=\mathbb{R}J$ for a nonzero $J \in \og(q)$, then the isomorphism type of $\n$ is determined by a single
number $d = \tfrac12 \rk J >0$.

To prove that $\n$ is an Einstein nilradical, it suffices to produce an Einstein derivation $\Phi$ and an inner product
$\<\cdot, \cdot\>$ on $\n$ (the nilsoliton inner product) in such a way that \eqref{eq:ricn} is satisfied for some number
$c<0$. Choose a basis $\{X_i\}$ for $\b$ in such a way that the matrix of $J$ has the form
$J =\left(\begin{smallmatrix} J_{2d} & 0 \\0 & 0 \end{smallmatrix}\right)$, where
$J_{2d} =\left(\begin{smallmatrix} 0 & I_d \\-I_d & 0 \end{smallmatrix}\right)$.
Introduce the inner product on $\b$ in such a way that the basis $\{X_i\}$ is orthonormal. Denote
$\b_1 = \Span(X_1, \ldots, X_{2d}), \; \b_2 = \b_1^{\perp} = \Ker J$ and  $l = q - 2d = \dim \b_2$.
By the result of \cite[Proposition 2.9]{GK}, we can assume that $l>0$ .

Let $F$ be an endomorphism of $\b$ whose eigenspaces are $\b_1$ and $\b_2$, with the corresponding eigenvalues
$\mu_1 \ne \mu_2$, respectively. The operator $g(F)$ acting on $\og(q)$ according to \eqref{eq:der2step} is semisimple
and real, with the eigenvalues $2\mu_1, \mu_1+\mu_2, 2\mu_2$ whose corresponding eigenspaces are
\begin{equation*}
    L_{2\mu_1} = \Bigl\{\left(\begin{matrix} K & 0 \\0 & 0 \end{matrix}\right), \; K \in \og(2d)\Bigr\}, \;
    L_{\mu_1+\mu_2} = \Bigl\{\left(\begin{matrix} 0 & T \\-T^t & 0 \end{matrix}\right) \Bigr\}, \;
    L_{2\mu_2} = \Bigl\{\left(\begin{matrix} 0 & 0 \\0 & N \end{matrix}\right) , \; N \in \og(l)\Bigr\}.
\end{equation*}
Clearly, $J \in L_{2\mu_1}$, so $W^\perp = L_{2\mu_1}' \oplus L_{\mu_1+\mu_2} \oplus L_{2\mu_2}$, where
$L_{2\mu_1}' = L_{2\mu_1} \cap J^\perp$ (recall that $\perp$ refers to the inner product in $\og(q)$; the above
decomposition of $W^\perp$ is orthogonal with respect to that inner
product).

Introduce an inner product on $\m$ in such a way that $\m$ admits an orthogonal decomposition
$\m = \m_{11} \oplus \m_{12} \oplus \m_{22}$ with the following property: for every $Z \in \m_{ij}$, the operator
$J_Z$ belongs to $L_{\mu_i+\mu_j}$. Choose an orthonormal basis $\{Z_1,\ldots, Z_{D-1}\}$ for $\m$ such that
$\m_{11} = \Span(Z_1, \ldots , Z_{d(2d-1)-1}),\, \m_{12} = \Span(Z_{d(2d-1)}, \ldots , Z_{d(2d-1)-1+2dl}),\,
\m_{22} = \Span(Z_{d(2d-1)-1+2dl}, \ldots , Z_{D-1})$. Equation \eqref{eq:ricn} gives
\begin{equation}\label{eq:p=D-1ricn}
\sum\nolimits_{k=1}^{D-1} J_k^2 = 2(c I_q + F), \qquad
\Tr (J_rJ_s^t) = 4 (c+\mu_i+\mu_j) \K_{rt}, \quad \text{for} \; Z_r \in \m_{ij},
\end{equation}
where $I_q$ is the identity matrix (we abbreviate $J_{Z_k}$ to $J_k$). Now choose any bases $\{J_r\}$ in
each of the subspaces $L_{2\mu_1}', L_{\mu_1+\mu_2}, L_{2\mu_2}$ such that the second equation of
\eqref{eq:p=D-1ricn} is satisfied. Then an easy computation shows that the first equation of \eqref{eq:p=D-1ricn} is
equivalent to
$$
2d(c+\mu_1+\mu_2) + (l-1) (c+2\mu_2) = -(c+\mu_2), \quad l(c+\mu_1+\mu_2) + (2d-1-1/d) (c+2\mu_1) = -(c+\mu_1).
$$
Solving this we find, up to scaling:
$$
\mu_1 = (q^2-q+2)d-2q, \quad \mu_2=(q^2-q+2)d-2q+1, \quad -c = (2q^2-3q+5)d-4q+2.
$$
To check that \eqref{eq:p=D-1ricn} is satisfied, it remains to show that $\mu_i < -c < \mu_i + \mu_j$.
As $d \ge 1$ and $q=2d+l > 2d$, this is indeed the case, with the only exception: when $q=3$ and
$d=1, \quad -c = 2\mu_1 = 4$. However, when $d=1, \; L_{2\mu_1}' = 0$, so it is sufficient to check that
$-c < \mu_1 + \mu_2$, which is true.
\end{proof}

\begin{remark}
In view of Lemma~\ref{l:p=D-1}, the fact that any two-step nilpotent Lie
algebra of type $(1, q)$ is an Einstein nilradical and
\cite[Proposition 2.9 (iv)]{GK}, one may wonder, whether a
two-step nilpotent Lie algebra dual to an Einstein nilradical is an Einstein nilradical by itself. This is
not true in general, as can be seen from the example of the two-step nilpotent Lie algebra $\n$ attached to the graph
$\mathcal{G}_{2,2,0}$ \cite[Figure 1]{LW}. By \cite[Proposition~5.6]{LW}, $\n$ is not an Einstein nilradical.
The dual algebra $\overline{\n}$ is attached to the graph which has the same vertex set as
$\mathcal{G}_{2,2,0}$ and whose edge set is the complement to that of $\mathcal{G}_{2,2,0}$. Applying
\cite[Theorem 5.3]{LW} one can easily see that $\overline{\n}$ is an Einstein nilradical.
\end{remark}

\begin{lemma}\label{l:small}
1. Any Einstein metric solvable Lie algebra whose nilradical is a two-step nilpotent Lie algebra of type
$(p,q), \; q \le 4$ is standard.

2. Any two-step nilpotent Lie algebra of type $(2,5)$ is an Einstein nilradical.
\end{lemma}
\begin{proof}
1. If $q \le 3$, then $\dim \n \le 6$ and the claim follows from \cite[Theorem 5.1]{L3} and \cite[Theorem 3.1]{W}.

Any two-step nilpotent Lie algebra of type $(p,4)$ with $p=1, 2, 5, 6$ is an Einstein nilradical (this follows from
\cite[Theorem 5.1]{L3}, \cite[Theorem 3.1]{W}, Lemma~\ref{l:p=D-1}, and \cite[Proposition 2.9(iii)]{GK}, respectively),
and the claim then follows from Theorem~\ref{th:nil}.

It remains to consider the cases $p=3, 4$.
Let $(\g, \<\cdot, \cdot\>)$ be an Einstein metric solvable Lie algebra whose nilradical $\n$ is two-step nilpotent
of type $(p,4),\; p = 3, 4$. Let $\n = \b \oplus \m$ be the orthogonal decomposition, with $\m = [\n, \n]$, and let
$\ag = \n^\perp$.

For every $Y \in \ag$, define the operators $F(Y): \b \to \b$ and $M(Y):\m \to \m$ by $F(Y) = \pi_{\b} \circ \ad_{Y|\b}$
and $M(Y) = \ad_{Y|\m}$, where $\pi_{\b}$ is the orthogonal projection to $\b$. Then both
$\{F(Y)\, |\, Y \in \ag\}$ and $\{M(Y)\, |\, Y \in \ag\}$ are commuting families of operators, and the same is true
for their semisimple parts $F(Y)^S, M(Y)^S$ (see Section~\ref{ss:2steppre}). Note also that $M(Y)$ is completely
determined by $F(Y)$ by \eqref{eq:der2step}. In particular, if $F(Y)^S$ is real, then $M(Y)^S$ is real, hence
$B(Y, Y) \ge 0$ (where $B$ is the Killing form of $\g$).

The linear space $F(Y)^S$ is a commuting family of semisimple endomorphisms of $\mathbb{R}^4 = \b$. We have a
direct sum decomposition $\b = \b_1 \oplus \ldots \oplus \b_m$ on invariant subspaces, such that every $\b_i$ is either
one- or two-dimensional. When $\dim \b_i = 1, \; F(Y)^S_{|\b_i} = l_i(Y)\, \id_{|\b_i}$ for some linear form
$l_i$ on $\ag$. When $\dim \b_i = 2$, we can choose a basis $X_i^1, \, X_i^2$ for $\b_i$ such that
$F(Y)^S_{|\b_i} = c_i(Y)\,\id_{|\b_i} + d_i(Y) J_i$ for some linear forms $c_i, d_i$ on $\ag$, with
$J_i: \b_i \to \b_i$ defined by $J_i(X_i^1) = X_i^2, \; J_i(X_i^2) = -X_i^1$.

If all the $\b_i$'s are of dimension one (all the eigenvalues of all the $F(Y)$'s are real), then
$B(Y, Y) \ge 0$, for all $Y \in \ag$. If all the $\b_i$'s except one are one-dimensional,
then $\ag$ contains a subspace of codimension one such that for every $Y$ from that subspace all the eigenvalues
of $F(Y)$ are real, so the index of $B$ is at most one. In the both cases, the claim follows from
\cite[Lemma 4.7]{H}. The same is true if we have two two-dimensional subspaces $\b_1,\, \b_2$ and the linear forms
$d_1, \; d_2$ are proportional. The only remaining case is therefore when
$\b = \b_1 \oplus \b_2, \; \dim \b_1 = \dim \b_2 = 2$ and the forms $d_1, \, d_2$ are linearly independent.
As for a generic $Y \in \ag$, all the (complex) eigenvalues of $F(Y)$ are distinct, all the
$F(Y)$'s are semisimple (and hence all the $M(Y)$'s are semisimple). Fix the basis
$X_1^1, \, X_1^2, \, X_2^1, \, X_2^2$ for $\b$. Then the matrix of $F(Y)$ has the form
$F(Y)  = \left(\begin{smallmatrix} c_1(Y) I_2 + d_1(Y) J & 0 \\0 & c_2(Y) I_2 + d_2(Y) J \end{smallmatrix}\right)$,
where $J = \left(\begin{smallmatrix} 0 & 1 \\ -1 & 0 \end{smallmatrix}\right)$. The operators
$g(F(Y)): \og(4) \to \og(4)$ defined by \eqref{eq:der2step} are commuting and semisimple
and have four eigenspaces:
\begin{gather*}
L_1 = \mathbb{R}\begin{pmatrix} J & 0 \\ 0 & 0 \end{pmatrix},\quad
L_2 = \mathbb{R}\begin{pmatrix} 0 & 0 \\ 0 & J \end{pmatrix},\\
L_3 = \Span\Bigl(\begin{pmatrix} 0 & I_2 \\ -I_2 & 0 \end{pmatrix},\,
\begin{pmatrix} 0 & J \\ J & 0 \end{pmatrix}\Bigr), \quad
L_4 = \Span\Bigl(\begin{pmatrix} 0 & I_{1,1} \\ -I_{1,1} & 0 \end{pmatrix},\,
\begin{pmatrix} 0 & I_{1,1}J \\ -I_{1,1}J & 0 \end{pmatrix}\Bigr),
\end{gather*}
where
$I_{1,1} = \left(\begin{smallmatrix} 1 & 0 \\ 0 & -1 \end{smallmatrix}\right)$, with the corresponding eigenvalues
$2c_1(Y)$, $2c_2(Y), \; c_1(Y) + c_2(Y) \pm i(d_2(Y) - d_1(Y))$, $c_1(Y) + c_2(Y) \pm i(d_2(Y) + d_1(Y))$
(the first two span a single two-dimensional eigenspace, when $c_1(Y) = c_2(Y)$).

Suppose $p=3$. If $c_1(Y) = c_2(Y)$, any
$g(F(\ag))$-invariant three-dimensional subspace of $\og(4)$ is a direct sum
of a one-dimensional subspace from $L \subset L_1 \oplus L_2$ and one of the spaces $L_3, L_4$. A direct check
shows that if $L \not\subset L_1 \cup L_2$, then $\n \in \Og(3,4)$, and the claim follows from Lemma~\ref{l:opqstandard}
(in fact, such an $\n$ is an Einstein nilradical). If $L = L_1$
or $L = L_2$, then $\n$ is an Einstein nilradical, as it is shown below.

If $c_1(Y) \ne c_2(Y)$, then every three-dimensional $g(F(\ag))$-invariant subspace of $\og(4)$ is a direct sum
of one of the $L_1, L_2$ and one of the $L_3, L_4$. Without loss of generality (permuting the basis vectors and changing
the sign of some of them, if necessary), we can take $W = L_1 \oplus L_3$. Then $\n$ is defined
by the commutator relations
$[X_1, X_2] = Z_1,\; [X_1, X_3] = [X_2, X_4] = Z_2,\; [X_1, X_4] = [X_3, X_2] = Z_3$ and is an Einstein nilradical
(the nilsoliton inner product is defined by taking the vectors $X_i, Z_k$ orthonormal).

Let now $p=4$. There are two $g(F(\ag))$-invariant four-dimensional subspaces of $\og(4)$: $L_1 \oplus L_2 \oplus L_3$
(choosing $L_4$ instead of $L_3$ gives an isomorphic algebra) and $L_3 \oplus L_4$.
Using \cite[Theorem 1]{P} we find that in the both cases, $\n$ is an Einstein nilradical: in the first case, $\n$ is
isomorphic to the Lie algebra given by the commutator relations
$[X_1, X_2] = \sqrt{2} Z_1,\; [X_3, X_4] = \sqrt{2} Z_2,\; [X_1, X_3] = [X_2, X_4] = Z_3$,
$[X_1, X_4] = [X_2, X_3] = Z_4$,
in the second one --- by the commutator relations
$[X_1, X_3] = Z_1$, $[X_1, X_4] = Z_2$, $[X_2, X_3] = Z_3,\; [X_2, X_4] = Z_4$
(the nilsoliton inner product in the both cases is defined by taking $X_i, Z_k$ orthonormal).
The claim then follows from Theorem~\ref{th:nil}.

2. A two-step nilpotent Lie algebra $\n$ of type $(2,5)$ is determined by a point $W \in G(2,\og(5))$. Let
$W=\Span(J_1, J_2)$, with $J_1, J_2 \in \og(5)$ linearly independent. If $\Ker J_1 \cap \Ker J_2 \ne 0$, the algebra
$\n$ is decomposable: it is a direct sum a two-step nilpotent Lie algebra of a smaller dimension and an abelian ideal.
Any such $\n$ is an Einstein nilradical, as it follows from
\cite[Theorem 3.1]{W}, \cite[Theorem 5.1, Proposition 3.3]{L3}.

Assume $\Ker J_1 \cap \Ker J_2 = 0$. Then by \cite[Theorem 5.1]{LR}, $\n$ is isomorphic to one of the two algebras
given by the following $W=\Span(J_1, J_2)$:
\begin{equation*}
J_1 = \begin{pmatrix} J & 0 \\ 0 & 0_3 \end{pmatrix}, \;
J_2 = \begin{pmatrix} 0_1 & 0 & 0\\ 0 & J & 0 \\ 0 & 0 & J \end{pmatrix},
\quad \text{or} \quad
J_1 = \begin{pmatrix} 0 & J & 0 \\ J & 0 & 0 \\ 0 & 0 & 0_1 \end{pmatrix}, \;
J_2 = \begin{pmatrix} 0 & 0 & J\\ 0 & 0_1 & 0 \\ J & 0 & 0 \end{pmatrix}
\end{equation*}
where $J = \left(\begin{smallmatrix} 0 & 1 \\ -1 & 0 \end{smallmatrix}\right)$.

Using \cite[Theorem 1]{P} we find that in the both cases, $\n$ is an Einstein nilradical: in the first case, $\n$ is
isomorphic to the Lie algebra given by the commutator relations
$[X_1, X_2] = 2 Z_1,\; [X_1, X_5] = \sqrt{2} Z_2$, $[X_3, X_4] = 2 Z_2$,
in the second one --- by the commutator relations
$[X_1, X_4] = Z_1,\; [X_2, X_4] = Z_2$, $[X_2, X_3] = \sqrt{2} Z_1,\; [X_1, X_5] = \sqrt{2} Z_2$
(the nilsoliton inner product in the both cases is defined by taking $X_i, Z_k$ orthonormal).
\end{proof}

As a consequence of Lemma~\ref{l:small} we obtain that any Einstein metric solvable Lie algebra with a two-step
nilradical of dimension at most seven is standard.

\section{Low-dimensional nilradicals. Proof of Theorem~\ref{th:dim7}}
\label{s:dim7}

The second assertion of Theorem~\ref{th:dim7} follows trivially from the first one (and Lemma~\ref{l:nilradvsderived}).
Let $(\g, \<\cdot,\cdot\>)$ be an Einstein metric solvable Lie algebra with the nilradical $\n, \; \dim \n \le 7$.
If $\dim \n \le 6$, the claim follows from the fact that any such $\n$ is an Einstein nilradical
(\cite[Theorem 3.1]{W}, \cite[Theorem 5.1]{L3}) and Theorem~\ref{th:nil}. Suppose $\dim \n = 7$. If $\n$ is
decomposable (is a direct sum of nonzero ideals), then
each of the summands is an Einstein nilradical, hence such is $\n$ \cite[Theorem 4]{P}.
Denote $\m = [\n, \n],\; \dim \m = p$. Then $\n$ is generated by $7-p$ elements. If $p \ge 4$, the claim follows from
Theorem~\ref{th:filiform}. If $p=0$ or $p=1, \; \n$ is an Einstein nilradical: in the first case, $\n$ is abelian,
in the second one, $\n$ is a Heisenberg algebra (if it is indecomposable).
We have two remaining cases to consider: $p=2$ and $p=3$. In each of these
cases, if $\n$ is two-step nilpotent, the claim follows from Lemma~\ref{l:small}.

So it suffices to consider only those seven-dimensional nilpotent algebras, which are indecomposable and not
two-step nilpotent, and whose derived algebra has dimension $2$ or $3$.

First suppose $p=2$. Let $\m = \Span(Z_1, Z_2)$, with $Z_2$ in the center of $\n$. Then
$\z(\n) = \mathbb{R}Z_2$. Indeed, $Z_1$ is not in $\z(\n)$ (as otherwise $\n$ is two-step), and no vector from
$\n \setminus \m$ is in $\z(\n)$ (as otherwise the span of that vector is a direct summand, so $\n$ is decomposable). Let
$\b$ be a subspace of $\n$ complementary to $\m$.
Consider the operator
$\ad_{Z_1}$ (which is nonzero). Clearly, $\ad_{Z_1|\m}=0$ and  $\ad_{Z_1}(\b) \subset \m$. What is more,
$\ad_{Z_1}(\b) = \mathbb{R}Z_2$ (for otherwise $\n$ is not nilpotent), so we get a nonzero one-form $l$ on $\b$
such that $[Z_1, X] = l(X) Z_2$. Define the two-forms $\omega^1, \omega^2$ on $\b$ by
$[X, Y] = \omega^1(X,Y) Z_1 + \omega^2(X,Y) Z_2$. From the Jacobi identity, $\omega^1(X,Y) = 0$ for all
$X, Y \in \Ker l$, hence, $\omega^1 = l \wedge l_1$ for some $l_1 \nparallel l$ (note that $\omega^1 \ne 0$, as
otherwise $Z_1 \notin \m$),
and in particular, $\rk \, \omega^1 = 2$. Moreover, the intersection $\Ker \omega^1 \cap \Ker \omega^2$ is the set of
generators whose bracket with any generator is zero, hence it lies in $\z(\n)$, so
$\Ker \omega^1 \cap \Ker \omega^2 = 0$. In particular, the forms $\omega^1$ and $\omega^2$ are not proportional. From
\cite[Theorem 5.1]{LR}, a pair of non-proportional $5 \times 5$ skew-symmetric matrices without a common kernel and such
that one of them has rank two can be reduced by a choice of the basis to the form
\begin{equation*}
J_1 = \begin{pmatrix} 0 & 1 & 0 & 0 \\ -1 & 0 & 0 & 0 \\ 0 & 0 & 0_1 & 0 \\ 0 & 0 & 0 & 0_2 \end{pmatrix}, \quad
J_2 = \begin{pmatrix} 0 & 0 & 1 & 0 \\ 0 & 0_1 & 0 & 0 \\ -1 & 0 & 0 & 0 \\ 0 & 0 & 0 & J \end{pmatrix},
\end{equation*}
where $J = \left(\begin{smallmatrix} 0 & 1 \\ -1 & 0 \end{smallmatrix}\right)$. Let $X_1, \ldots, X_5$ be the
corresponding basis for $\b$. As $\omega^1$ vanishes on any two vectors from $\Ker l$, we obtain
$l(X_i) = 0, \, i = 3,4,5$. Then $\n$ is defined by the following commutator relations:
$$
[X_1, X_2] = Z_1, \; [X_1, X_3] = Z_2, \; [X_4, X_5] = Z_2, \; [X_1, Z_1] = a Z_2, \; [X_2, Z_1] = b Z_2,
$$
where $a$ and $b$ are not simultaneously zeros. If $b = 0$, then $a \ne 0$, so the vector $X_3 - Z_1/a$ lies in $\z(\n)$,
which is a contradiction. If $b \ne 0$, we can replace $X_1$ by $X_1 - (a/b)X_2$, which gives after scaling the following
set of relations for $\n$: $[X_1, X_2] = Z_1, \; [X_1, X_3] = Z_2, \; [X_4, X_5] = Z_2, \; [X_2, Z_1] =  Z_2$.
By \cite[Theorem~1]{P}, we obtain that such an $\n$ is an Einstein nilradical, as it is isomorphic to the (metric) Lie
algebra given by the commutator relations
$[X_1, X_2] = 2 Z_1, \; [X_1, X_3] = Z_2, \; [X_4, X_5] = \sqrt{3} Z_2, \; [X_2, Z_1] = \sqrt{3} Z_2$,
with the nilsoliton inner product defined by taking $X_i, Z_k$ orthonormal. The claim now follows from
Theorem~\ref{th:nil}.

Now suppose $p=3$. We combine the approaches used in the proof of Theorem~\ref{th:split},
Theorem~\ref{th:filiform}, and Lemma~\ref{l:small}. As it was shown in the proof of Theorem~\ref{th:split},
there exists a subspace $\ag \subset \g$ complementary to $\n$ such that for any
$Y_1, Y_2 \in \ag, \; \ad_{Y_1}^SY_2 = [\ad_{Y_1}^S, \ad_{Y_2}^S] = 0$, and the map $Y \to \ad_Y^S$ is linear on $\ag$
(such an $\ag$ is the complement, in a Cartan subalgebra $\h$ of $\g$, to $\h \cap \n$).

The set of operators $(\ad_{Y}^S)_{|\n}$ is a torus, a commuting family of semisimple derivations of $\n$. As
$\m = [\n, \n]$ is a characteristic ideal, it is complemented in $\n$ by a linear subspace $\b$ invariant with
respect to all of the $\ad_{Y}^S$.

As in the proof of Theorem~\ref{th:filiform}, define a linear map $\T: \ag \to \End(\b)$ by
$\T(Y) = \pi_\b \circ (\ad_Y^S)_{|\n}$ (where $\pi_\b$ is the linear projection on $\b$). Then $\T(\ag)$ is a commuting
family
of semisimple endomorphisms of $\b$. If all the eigenvalues of $\T(Y)$ are real, then all the eigenvalues of $\T(Y)$ are
real, and so $B(Y,Y) = \Tr (\ad_Y^S)^2 \ge 0$ (where $B$ is the Killing form of $\g$). By \cite[Lemma 4.7]{H}, the metric
algebra $(\g, \<\cdot,\cdot\>)$ is standard, provided $\ag$ contains a subspace $\ag'$ of codimension at most one such
that $B(Y,Y) \ge 0$ for all $Y \in \ag'$. In particular, this will be the case, if $\ag$ contains a subspace $\ag'$ of
codimension at most one such that all the eigenvalues of $\T(Y)$ are real for all $Y \in \ag'$.
The only possible way to avoid that for a commuting family of semisimple endomorphisms of the four-dimensional space
$\b$ is when $\b$ splits into a direct sum of two two-dimensional subspaces $\b_1, \b_2$ invariant with respect to
$\T(Y)$ and corresponding to the complex eigenvalues $c_j(Y) \pm i d_j(Y), j=1,2$, with the linear forms $d_1, d_2$ on
$\ag$ linearly independent (see the proof of assertion 1 of Lemma~\ref{l:small}). Choose the bases $X_j^1, X_j^2$
for $\b_j$ such that $X_j^1 + i X_j^2$ is an eigenvector of $\T(Y)$ with the eigenvalue $c_j(Y) + i d_j(Y)$.

Then the weights for the torus $\{\ad_{Y|\n}^S\}$ are of the form $k_1 c_1(Y)+ k_2c_2(Y) + i(l_1d_1(Y) + l_2d_2(Y))$ for
some integers $k_1, k_2, l_1, l_2$. As $d_1$ and $d_2$ are linearly independent, two such weights can only be equal
if their $l_1$'s are equal and $l_2$'s are equal. The corresponding weight space has an even dimension, unless
$l_1 = l_2 = 0$.
As $\n$ is assumed to be not two-step nilpotent and is generated by $\b$, the space $[[\b,\b],\b]$ must be nonzero.
Note that $[[\b_1, \b_2], \b_1] = 0$, as otherwise the three-dimensional space $\m$ would contain two
even-dimensional weight spaces corresponding to the different (and not complex conjugate) weights. Similarly,
$[[\b_1, \b_2], \b_2] = 0$, and by the Jacobi identity, $[[\b_1, \b_1], \b_2] = [[\b_2, \b_2], \b_1] = 0$. The only
three-folded brackets which could be nonzero are therefore $[[\b_1, \b_1], \b_1]$, or $[[\b_2, \b_2], \b_2]$. Assume
that say $[[\b_1, \b_1], \b_1] \ne 0$. Then the vector $Z = [X_1^1, X_1^2]$ is nonzero, and $\m$ is a direct sum of the
one-dimensional weight space $\mathbb{R}Z$, with the weight $2c_1(Y)$, and the two-dimensional weight space
$[Z, \b_1]$, with the weight $3c_1(Y) \pm d_1(Y)$. As no subspace from $[\b_1, \b_2]$ can have such weights, we obtain
that $[\b_1, \b_2]=0$. Since $\b$ generates $\n$, this implies that $\n$ is decomposable: it is a direct sum of the
subalgebras (in fact, the ideals) generated by the $\b_j$'s, which is a contradiction.


\end{document}